\RequirePackage{fix-cm}
\documentclass[smallextended]{svjour3}       % onecolumn (second format)
\smartqed  % flush right qed marks, e.g. at end of proof

\usepackage{todonotes}
\usepackage{graphicx}
\usepackage{amssymb}
\usepackage{amsmath}
\usepackage{url}
\usepackage{dsfont}
\usepackage{ntheorem}
\makeatletter
\let\cl@chapter\undefined
\makeatletter
\usepackage[capitalise]{cleveref}
\usepackage{enumerate}
\usepackage{color}

\newtheorem{assumption}{Assumption}
\newtheorem{thm}{Theorem}
\newtheorem{lem}[thm]{Lemma} 
\newtheorem{prop}[thm]{Proposition} 
\newtheorem{defi}[thm]{Definition} 
\newtheorem{cor}[thm]{Corollary}
\newtheorem{examp}[thm]{Example}
\newtheorem{rem}[thm]{Remark}
\DeclareMathOperator*{\argmin}{argmin}
\DeclareMathOperator*{\sgn}{sgn}
\newcommand{\inv}{^{-1}}             %inversion
\newcommand{\gobs}{g^\mathrm{obs}}

\newcommand{\nat}{\mathbb{N}_0}

\newcommand{\Yspace}{\mathbb{Y}}
\newcommand{\smax}{s_\textrm{max}}

\newcommand{\xh}{\hat{x}_\alpha}
\newcommand{\new}[1]{#1}
\newcommand{\old}[1]{}
\newcommand{\lspace}[2]{\ell_{#1}^{#2}}
\newcommand{\lnorm}[3]{\left\| {#1} \right\|_{{#2},{#3}}}
\newcommand{\Bn}[4]{\| {#1} \|_ {B^{#2}_{{#3},{#4}}} }
\newcommand{\bn}[4]{\| {#1} \|_{{#2},{#3},{#4}}}
\newcommand{\wav}{\mathcal{S}}
\newcommand{\bspace}[3]{b^{#1}_{{#2},{#3}}}
\newcommand{\Bspace}[3]{B^{#1}_{{#2},{#3}}(\Omega)}
\newcommand{\wei}[1]{ {\underline{#1}} }
\begin{document}

\title{Maximal Spaces for Approximation Rates in $\ell^1$-regularization%\thanks{Grants or other notes
}


\author{Philip Miller          \and
       Thorsten Hohage %etc.
}


\institute{Institute for Numerical and Applied Mathematics, University of Göttingen, Germany\\
              \email{p.miller@math.uni-goettingen.de}           %  \\
}

\date{Received: date / Accepted: date}

\maketitle

\begin{abstract}

We study Tikhonov regularization for possibly nonlinear inverse problems with weighted $\ell^1$-penalization. 
The forward operator, mapping from a sequence space to an arbitrary Banach space, typically an $L^2$-space, is assumed to 
satisfy a two-sided Lipschitz condition with respect to a weighted $\ell^2$-norm 
and the norm of the image space.  
We show that in this setting approximation rates of arbitrarily high H\"older-type order in the regularization parameter 
can be achieved, and we characterize maximal subspaces of sequences on which these rates are attained.  
On these subspaces the method also converges with optimal rates in terms of the noise level with the discrepancy principle 
as parameter choice rule. 
Our analysis includes the case that the penalty term is not finite at the exact solution ('oversmoothing'). 
As a standard example we discuss wavelet regularization in Besov spaces $B^r_{1,1}$. 
\new{In this setting we demonstrate in numerical simulations for a parameter identification problem in a differential equation 
that our theoretical results correctly predict improved rates of convergence for piecewise smooth unknown coefficients.}  

\end{abstract}

\section{Introduction} 
In this paper we analyze numerical solutions of ill-posed operator equations 
\[
F(x)=g
\] 
with a (possibly nonlinear) forward operator $F$ mapping sequences 
$x=(x_j)_{j\in\Lambda}$ indexed by a countable set $\Lambda$ to a Banach space $\Yspace$. 
We assume that only indirect, noisy observations $\gobs\in \Yspace$ of the unknown solution $x^\dagger\in\mathbb{R}^\Lambda$ 
are available satisfying a deterministic error bound ${\|\gobs-F(x^\dagger)\|_{\Yspace}\leq \delta}$. \\
For a fixed sequence of positive weights $(\wei r_j)_{j\in\Lambda}$ and a regularization parameter $\alpha>0$ 
we consider Tikhonov regularization of the form
\begin{align}\label{eq:Tik}
\xh \in \argmin_{x\in D} \left[ \frac{1}{2} \|\gobs-  F(x) \|_\Yspace^2 +\alpha \sum_{j\in \Lambda} \wei r_j |x_j|  \right]
\end{align}
where $D\subset \mathbb{R}^\Lambda$ denotes the domain of $F$. 
Usually, $x^\dagger$ is a sequence of coefficients with respect to some Riesz basis. 
One of the reasons why such schemes have become popular is that  the penalty term 
$\alpha \sum_{j\in \Lambda} \wei r_j |x_j|$ promotes sparsity of the estimators $\xh$ in the sense that only 
a finite number of coefficients of $\xh$ are non-zero. \new{ The latter holds true if $(\wei r_j)_{j\in\Lambda}$ decays not too fast relative to the ill-posedness of $F$ (see  \Cref{prop:existence} below). In contrast to \cite{lorenz:08} and related works, we do not require that $(\wei r_j)_{j\in\Lambda}$ is uniformly bounded away from zero. In particular, this allows us to consider Besov $B^0_{1,1}$-norm penalties given by wavelet coefficients.} \\
For an overview on the use 
of this method for a variety linear and nonlinear inverse problems in different fields of applications we 
refer to the survey paper \cite{JM:12} and to the special issue \cite{JMS:17}. 

\medskip
\paragraph*{\new{Main contributions:$\,$} }
The focus of this paper is on error bounds, i.e.\ rates of convergence of $\xh$ to $x^\dagger$ in some norm as the noise level $\delta$  
tends to $0$. \new{Although most results of this paper are formulated for general operators on weighted $\ell^1$-spaces, 
we are mostly interested in the case that $x_j$ are wavelet coefficients, and 
\begin{align}\label{eq:F_factorization}
F= G\circ \wav
\end{align}
is the composition of a corresponding wavelet synthesis operator $\mathcal{S}$ and an operator $G$ defined on a function space. 
We will assume that $G$ is \emph{finitely smoothing} in the sense that it satisfies a two-sided Lipschitz condition with respect 
to function spaces the smoothness index of which differs by a constant $a>0$ (see Assumption \ref{ass_operator_Besov} below and 
Assumption \ref{ass_operator} for a corresponding condition on $F$). The class of operators satisfying this condition includes 
in particular the Radon transform and nonlinear parameter identification problems for partial differential equations with 
distributed measurements. In this setting Besov $B^{r}_{1,1}$-norms can be written in the form of the penalty term in \eqref{eq:Tik}. 
In a previous paper \cite{HM:19} we have already addressed  sparsity promoting penalties in the form of Besov $B^0_{p,1}$-norms 
with $p\in [1,2]$. For $p>1$ only group sparsity in the levels is enforced, but not sparsity of the wavelet coefficients within each level. 
As a main result of this paper we demonstrate that the analysis in \cite{HM:19} as well as other works to be discussed below 
do not capture the full potential of estimators \eqref{eq:Tik}, i.e.\ the most commonly used case $p=1$: 
Even though the error bounds in \cite{HM:19} are optimal in a minimax sense, 
more precisely in a worst case scenario in $B^s_{p,\infty}$-balls, we will derive faster rates of convergence for an important 
class of functions, which includes piecewise smooth functions. The crucial point is that such functions 
also belong to Besov spaces with larger smoothness index $s$, but smaller integrability index $p<1$. 
%These theoretically predicted improved rates of convergence for such functions are indeed observed in our numerical experiments for a parameter identification problem in an elliptic differential equation. 
These results confirm the intuition that estimators of the form \eqref{eq:Tik}, which enforce sparsity also within each wavelet level, should perform well for signals which allow accuratele approximations by sparse wavelet expansions.

Furthermore, we prove a converse result, i.e.\ we characterize the maximal sets on which 
the estimators \eqref{eq:Tik} achieve a given approximation rate. These maximal sets turn out to be weak weighted $\ell^t$-sequences spaces or real interpolation spaces of Besov spaces, respectively. 
%In this sense, no further progress in this direction is possible.

Finally, we also treat the \emph{oversmoothing case} that $\sum_{j\in \Lambda} \wei r_j |x_j^\dagger|=\infty$, i.e.\ that the 
penalty term enforces the estimators $\xh$ to be smoother than the exact solution $x^\dagger$. 
For wavelet $B^r_{1,1}$ Besov norm penalties, this case may be rather unlikely for $r=0$, except maybe for delta peaks. 
However, in case of the Radon transform, our theory requires us to choose $r>\frac 12$, and more generally, mildly ill-posed problems 
in higher spatial dimensions require larger values of $r$ (see eq.~\eqref{eq:r_condition} below for details). 
Then it becomes much more likely that the penalty term fails to be finite 
at the exact solution, and it is desirable to derive error bounds also for this situation.    
So far, however, this case has only rarely been considered in variational regularization theory. 
}

\medskip
\paragraph*{\new{Previous works on the convergence analysis of \eqref{eq:Tik}:$\,$}} 
In the seminal paper \cite{DDD:04}   Daubechies, Defrise \& De Mol 
established the regularizing property of estimators of the form \eqref{eq:Tik} and suggested the so-called iterative thresholding algorithm 
to compute them. 
%As most of the literature on this topic, \cite{DDD:04} only treats the case of linear forward operators $F$.  
Concerning error bounds, 
the most favorable case is that the true solution $x^\dagger$ is sparse. In this case the convergence rate is linear in 
the noise level $\delta$, and sparsity of $x^\dagger$ is not only sufficient but (under mild additional assumptions) 
even necessary for a linear 
convergence rate (\cite{GHS:11}). However, usually it is more realistic to assume that $x^\dagger$ is only \emph{approximately} 
sparse in the sense that it can be well approximated by sparse vectors. 
%Two further special rates of convergence, 
%$\mathcal{O}(\delta^{1/2})$ and $\mathcal{O}(\delta^{2/3})$ have also been studied early since they are implied by simple 
%range conditions (see \cite{}). 
More general rates of convergence for linear operators $F$ were derived in \cite{burger2013convergence} 
based on variational source conditions. The rates were characterized in terms of the growth of the norms of the 
preimages of the unit vectors under $F^*$ (or relaxations) and the decay of $x^\dagger$. Relaxations of the first condition  
were studied in \cite{FH:15,FHV:15,flemming2018injectivity}. 
For error bounds in the Bregman divergence with respect to the $\ell^1$-norm we refer to \cite{BHK:18}. 
\old{In a previous work \cite{HM:19} we studied regularization by 1-homogeneous wavelet-Besov norms 
$B^0_{p,1}$ and derived order optimal convergence rates under bounds on the $B^{s}_{p,\infty}$-Besov norm of the true solution based 
on variational source conditions.  
Moreover, for $p=2$ we showed that boundedness of the $B^s_{2,\infty}$-norm is not only sufficient, but even necessary for certain 
approximation rates. 

For $p<2$, and in particular for the case $p=1$ corresponding to \eqref{eq:Tik}, we were not able to prove the necessity of our 
assumptions. This was the starting point of the present work: We asked for maximal sets of solutions $x^\dagger$ for 
which the estimators \eqref{eq:Tik} converge at a given H\"older rate.} 
In the context of statistical regression by wavelet shrinkage \old{such }maximal sets of signals for which a certain rate of convergence is achieved have been studied in detail (see \cite{CVKP:03}). 
\old{In contrast to \cite{HM:19}, the penalty term is imposed in sequence space rather than in function space, and our analysis is carried 
out in an abstract sequence space setting. 
Only at the end of this paper we introduce wavelet Besov norms and discuss the case that 
$\sum_{j\in \Lambda} \wei r_j |x_j| $ describes a $B^r_{1,1}$-norm.}  

\old{The basic assumption of our analysis is a two-sided Lipschitz condition for $F$ with respect to a weighted $\ell^2$-norm. 
For linear operators this is equivalent to the fact the $F$ is bounded and boundedly invertible from the weighted $\ell^2$-space 
to its image in $\Yspace$. As a main result we show that approximation rates can be equivalently characterized either 
by weak sequence spaces of solution or by variational source conditions. The latter conditions lead to convergence rates 
for noisy data as the noise level tends to $0$.}
 
%\new{In practical situations, especially in higher spatial dimensions, 
%especially when approximating functions defined on a domain of large dimension, 
%it may be an advantageous to use a strong penalty term. Here this means to take a large weights $\wei r_j$. At the same time most often the true solutions regularity is not known such that it is not unlikely that the true solution fails to admit a finite value in the penalty term. This motivates studying oversmoothing as mentioned above.} \\
\new{In the oversmoothing case} one difficulty is that neither variational source conditions nor source conditions based on the range of the 
adjoint operator are applicable. 
%The special case of diagonal operators in $\ell^1$-regularization has been discussed in \cite{GH:19}. 
Whereas oversmoothing in Hilbert scales has been analyzed in 
numerous papers (see, e.g., \cite{HM:18,HP:19,natterer:84}), the literature on oversmoothing for more general variational regularization 
is sparse. 
%One difficulty is that neither variational source conditions nor source conditions based on the range of the 
%adjoint operator are applicable. 
The special case of diagonal operators in $\ell^1$-regularization has been discussed in \cite{GH:19}. 
\new{In a very recent work, Chen, Hofmann \& Yousept \cite{CHY:21} have studied oversmoothing for finitely smoothing operators  
in scales of Banach spaces generated by sectorial operators.}

\medskip
\paragraph*{\new{Plan of the remainder of this paper:$\,$}}
\new{ 
In the following section we introduce our setting and assumptions and discuss two examples for which these assumptions 
are satisfied in the wavelet--Besov space setting \eqref{eq:F_factorization}. Sections \ref{sec:weaksequencespace}--\ref{sec:oversmoothing} 
deal with a general sequence space setting.} 
In Section \ref{sec:weaksequencespace} 
we introduce a scale of weak sequence spaces which can be characterized by the approximation properties of some hard thresholding 
operator. These weak sequence spaces turn out to be the maximal sets of solutions on which the method \eqref{eq:Tik} attains certain 
H\"older-type approximation rates. This is shown for the non-oversmoothing case in Section \ref{sec:non_oversmoothing} 
%using variational source conditions 
and for the oversmoothing case in Section \ref{sec:oversmoothing}. 
\new{In Section \ref{sec:Besov} we interpret our results in the previous sections in the Besov space setting, before we 
discuss numerical simulations confirming the predicted convergence rates in Section \ref{sec:num}.}
%we introduce wavelet Besov penalties of type $B^r_{1,1}$ and discuss 
%the consequences of our general results in this particularly relevant setting. 

\section{Setting, Assumptions, and Examples}\label{sec:setting}

\new{In the following we describe our setting in detail including assumptions which are used in many of the following results. 
None of these assumptions is to be understood as a standing assumption, but each assumption is referenced whenever it is needed.

\subsection{Motivating example: regularization by wavelet Besov norms}
In this subsection, which may be skipped in first reading, we provide more details on the motivating example \eqref{eq:F_factorization}:  
Suppose the operator $F$ is the 
composition of a forward operator $G$ mapping functions on a domain $\Omega$ to elements of the Hilbert space $\Yspace$ 
and a wavelet synthesis operator $\wav$. 
We assume that $\Omega$ is either a bounded Lipschitz domain in $\mathbb{R}^d$ or the $d$-dimensional torus $(\mathbb{R}/\mathbb{Z})^d$, 
and that we have a system $(\phi_{j.k})_{(j,k)\in\Lambda}$ of real-valued wavelet functions on $\Omega$. 
Here the index set $\Lambda:= \{(j,k) \colon j\in\nat, k\in\Lambda_j\}$ is composed of 
a family of finite sets $(\Lambda_j)_{j\in\nat}$ corresponding to levels $j\in \nat$, and the growths of the cardinality of these sets 
is described by the inequalities 
$2^{jd}\leq |\Lambda_j|\leq C_\Lambda 2^{jd}$ for some constant $C_\Lambda\geq 1$ and all $j\in\nat$. 

For $p,q \in (0,\infty)$ and $s\in \mathbb{R}$ we introduce sequence spaces 
\begin{align}\label{eq:besov_seq}
\begin{aligned}
\bspace spq &:=\left\{x\in\mathbb{R}^\Lambda \colon \bn xspq<\infty\right\}\qquad\mbox{  with }\\
\bn xspq^q &:= \sum_{j\in \nat} 2^{jq(s+\frac{d}{2}-\frac{d}{p})} \left( \sum_{k\in \Lambda_j} |x_{j,k}|^p \right)^\frac{q}{p}.  
\end{aligned}
\end{align}
with the usual replacements for $p=\infty$ or $q = \infty$. 
%In order to infer the sequence space assumption \ref{ass_operator} from \Cref{ass_operator_Besov} 
It is easy to see that $\bspace spq$ are Banach spaces if $p,q\geq 1$. Otherwise, 
if $p\in(0,1)$ or $q\in (0,1)$, they are quasi-Banach spaces, 
%the sets $\lspace {\wei \omega}  p$ are quasi-Banach spaces with quasi-norms
i.e.\ they satisfy all properties of a Banach space except for the triangle inequality, 
which only holds true in the weaker form $\lnorm {x+y} {\wei \omega}  p \leq C(\lnorm x {\wei \omega}  p+\lnorm y {\wei \omega}  p)$ with some $C>1$. %, here the best possible constant is $C=2^{1/p-1}$. 
We need the following assumption on the relation of the Besov sequence spaces to a family of Besov function spaces 
$\Bspace spq$ via the wavelet synthesis operator $(\wav x)(\mathbf{r}) := \sum_{(j,k)\in \Lambda} x_{j,k} \phi_{j,k}(\mathbf{r})$. 
\begin{assumption}\label{ass:wavelet}
Let $s_\textrm{max}>0$. 
Suppose that $(\phi_{j.k})_{(j,k)\in\Lambda}$ is a family of real-valued functions on $\Omega$ such that the synthesis operator  
\[ \wav \colon\bspace spq \rightarrow \Bspace spq \quad\text{ given by } x\mapsto \sum_{(j,k)\in \Lambda} x_{j,k} \phi_{j,k} 
\] 
is a norm isomorphism for all $s \in (-s_\textrm{max}, s_\textrm{max})$ and $p,q \in (0,\infty]$
satisfying ${s \in (\sigma_p-s_\textrm{max}, s_\textrm{max})}$ with $\sigma_p=\max\left\{d\left(\frac{1}{p}-1\right), 0 \right\}$. 
\end{assumption}
Note that $p\geq 1$ implies $\sigma_p=0$, and therefore $\mathcal{S}$ is a quasi-norm isomorphism for $|s|\leq s_\textrm{max}$ in this case. \\
We refer to the monograph \cite{triebel:08} for the definition of Besov spaces $\Bspace spq$, different types of Besov spaces 
on domains with boundaries, and the verification of \Cref{ass:wavelet}.

%We formulate the assumptions on the forward operator in terms of Besov function spaces $\Bspace spq$. 
%We refer, e.g., to \cite{BL:76,triebel:08} for details on these function spaces.
%The norm on $\Bspace spq$ will be denoted by $\Bn \cdot spq$.
As main assumption on the forward operator $G$ in function space we suppose that it is finitely smoothing in the following sense:
\begin{assumption} \label{ass_operator_Besov} 
Let $a>0$, $D_G\subseteq \Bspace {-a}22$ be non-empty and closed, $\Yspace$ a Banach space and $G \colon D_G \rightarrow \Yspace$ a map. Assume that there exists a constant $L\geq 1$ with 
\begin{align*}
\frac{1}{L}  \Bn {f_1-f_2}{-a}22\leq \|G(f_1)-G(f_2)\|_{\Yspace} \leq L \Bn {f_1-f_2}{-a}22  \quad\text{for all } f_1,f_2 \in D_G.
\end{align*}
\end{assumption} 
Recall that $\Bspace {-a}22$ coincides with the Sobolev space $H^{-a}(\Omega)$ with equivalent norms.  The first of these inequalities is 
violated for infinitely smoothing forward operators such as for the backward heat equation or for electrical impedance tomography. \\
In the setting of Assumptions \ref{ass:wavelet} and \ref{ass_operator_Besov} and for some fixed $r\geq 0$
we study the following estimators
\begin{align}\label{eq:Tik_besov}
\hat{f}_\alpha:= \wav \xh \quad\text{with}\quad
\xh \in \argmin_{x\in \wav\inv(D_G)} \left[ \frac{1}{2} \|\gobs-  G(\wav x) \|_\Yspace^2 +\alpha \bn x r11 \right].
\end{align}
 We recall two examples of forward operators satisfying \Cref{ass_operator_Besov} from \cite{HM:19} where further examples 
are discussed. } 
\begin{examp}[Radon transform]\label{ex:radon}
\new{
Let $\Omega\subset\mathbb{R}^d$, $d\geq 2$ be a bounded domain and $\Yspace = L^2(S^{d-1}\times \mathbb{R})$ with the unit sphere 
$S^{d-1}:=\{x\in \mathbb{R}^d:|x|_2=1\}$. The Radon transform, which occurs in computed tomography (CT) and positron emission tomography (PET), among others, is defined by 
\[
(Rf)(\theta,t):=\int_{\{x:x\cdot\theta=g\}} f(x)\,\mathrm{d}x,\qquad \theta\in S^{d-1},\, t\in\mathbb{R}.
\]
It satisfies \Cref{ass_operator_Besov} with $a=\frac{d-1}{2}$. }
\end{examp}
\begin{examp}[identification of a reaction coefficient]\label{ex:identification}
\new{
Let $\Omega\subset \mathbb{R}^d$, $d\in\{1,2,3\}$ be a bounded Lipschitz domain, 
and let $f:\Omega\to [0,\infty)$ and $g:\partial\Omega\to (0,\infty)$ be smooth functions. For $c\in L^{\infty}(\Omega)$ satisfying 
$c\geq 0$ we define the forward operator $G(c):=u$  by the solution of the elliptic boundary value problem 

\begin{align}\label{eq:bvp}
\begin{aligned}
&-\Delta u+ cu =f&& \mbox{in }\Omega,\\
&u=g &&\mbox{on }\partial \Omega.
\end{aligned}
\end{align}
Then \Cref{ass_operator_Besov} with $a=2$ holds true in some $L^2$-neighborhood of a reference solution $c_0\in L^{\infty}(\Omega)$, $c_0\geq 0$. 
(Note that for coefficients $c$ with arbitrary negative values uniqueness in the boundary value problem \eqref{eq:bvp} may fail 
and every $L^2$-ball contains functions with negative values on a set of positive measure, well-posedness of \eqref{eq:bvp} can still 
be established for all $c$ in a sufficiently small $L^2$-ball centered at $c_0$. This can be achieved by Banach's fixed 
point theorem applied to $u = u_0+(-\Delta + c_0)^{-1}(u(c_0-c))$ where $u_0:=G(c_0)$ and  
$(-\Delta + c_0)^{-1}\tilde{f}$ solves \eqref{eq:bvp} with $c=c_0$, $f=\tilde{f}$ and $g=0$, using the fact that 
$(-\Delta + c_0)^{-1}$ maps boundedly from $L^1(\Omega)\subset H^{-2}(\Omega)$ to $L^2(\Omega)$ for $d\leq 3$.) }
\end{examp}

\new{\subsection{General sequence spaces setting}}\label{sec:ass_sequence}
Let $p\in (0,\infty)$, and let $\wei \omega=(\wei \omega_j)_{j\in \Lambda}$ be a sequence of positive reals indexed 
by some countable set $\Lambda$. 
We consider weighted sequence spaces $\lspace {\wei \omega} p$ defined by 
\begin{align}\label{eq:defi_lpw}
 \lspace {\wei \omega}  p := \left\{ x\in \mathbb{R}^\Lambda \colon \lnorm x {\wei \omega}  p < \infty \right\} \quad\text{with } \quad \lnorm x {\wei \omega}  p:= \left(\sum_{j\in \Lambda} {\wei \omega}_j^p | x_j |^p \right)^\frac{1}{p}. 
\end{align} 
\new{Note that the Besov sequence spaces $\bspace spq$ defined in \eqref{eq:besov_seq}  are of this form 
if $p=q <\infty$, more precisely $\bspace spp =\lspace {\wei \omega_{s,p}}p$ with equal norm for $(\wei \omega_{s,p})_{(j,k)}= 2^{j(s+\frac{d}{2}-\frac{d}{p})}$.
Moreover,} the penalty term in  is given by $\alpha\lnorm \cdot {\wei r} 1$ with the sequence of weights ${\wei r} =(\wei r_j)_{j\in \Lambda}$. \new{Therefore, we obtain the penalty terms $\alpha\bn {\cdot}s11$ in  \eqref{eq:Tik_besov} for 
the choice $\wei r_{j,k} := 2^{j(r-\frac{d}{2})}$.}

We formulate a two-sided Lipschitz condition for forward operators $F$ on general sequence spaces and argue that it follows from Assumptions \ref{ass:wavelet} and \ref{ass_operator_Besov} in the Besov space setting.
\begin{assumption} \label{ass_operator}  
$\wei a=(\wei a_j)_{j\in \Lambda}$ is a sequence of positive real numbers with \linebreak ${\wei a_j \wei r_j\inv \rightarrow 0}$.\footnote{This notion means that for every $\varepsilon>0$ all but finitely many $j\in \Lambda$ satisfy  $\wei a_j \wei r_j\inv\leq \varepsilon$.} 
Moreover, $D_F\subseteq \lspace  {\wei a} 2 $ is closed with $D_F\cap \lspace {\wei r} 1 \neq \emptyset$ and  
there exists a constant $L>0$ with 
\[ \frac{1}{L} \lnorm {x^{(1)} - x^{(2)}}{\wei a}2 \leq \|F(x^{(1)})-F(x^{(2)}) \|_\Yspace 
\leq L  \lnorm {x^{(1)} - x^{(2)}}{\wei a}2  \] 
for all $x^{(1)},x^{(2)}\in D_F$.
\end{assumption} 
\new{
Suppose Assumptions \ref{ass:wavelet} and \ref{ass_operator_Besov} hold true, and let 
\begin{subequations}\label{eqs:r_condition}
\begin{align}\label{eq:r_condition}
&\frac{d}{2}-r<a<\smax,\\
\label{eq:r_nonneg}
&r\geq 0,\\
&\wav\inv(D_G)\cap \bspace r11\neq \emptyset. 
\end{align} 
\end{subequations}
With $\wei a _{j,k}:=2^{-ja}$ and $\wei r_{j,k} := 2^{j(r-\frac{d}{2})}$
we have  $\lspace {\wei a}2= \bspace {-a}22$ and $\lspace {\wei r}1 = \bspace r11$. Then $\wei a _{j,k} \wei r_{j,k}\inv  \rightarrow 0$. As $\mathcal{S}\colon \bspace {-a}22 \rightarrow \Bspace {-a}22$ is a norm isomorphim $D_F:= \mathcal{S}\inv (D_G)$ is closed, and $F:=G\circ \mathcal{S}:D_F\to \Yspace$ satisfies the two-sided Lipschitz condition above. } 

In some of the results we also need the following assumption on the domain $D_F$ of the map $F$. 
\begin{assumption} \label{ass_domain}
$D_F$ is closed under coordinate shrinkage. That is 
$x\in D_F$ and $z\in \lspace {\wei a}2$ with $|z_j|\leq |x_j|$ \new{and $\sgn z_j\in\{0,\sgn x_j\}$}
for all $j\in \Lambda$ implies $z\in D_F. $
\end{assumption}
\new{
Obviously, \Cref{ass_domain} is satisfied if $D_F$ is a closed ball $\{x\in \lspace {\wei a} 2 : \lnorm x {\wei \omega} p \leq \rho\}$ in some 
$\lspace {\wei \omega} p$  space centered at the origin. 

Concerning the closedness condition in \Cref{ass_operator}, note 
that such balls are always closed in $\lspace {\wei a} 2$ as the following argument shows: 
Let $ x^{(k)}\to  x$ as $k\to\infty$ in $\lspace {\wei a} 2$ and $\lnorm { x^{(k)}} {\wei \omega} p\leq \rho$ for all $k$. 
Then $x^{(k)}$ converges pointwise to $ x$, and hence 
$\sum_{j\in\Gamma} \wei \omega_j^p |x_j|^p = \lim_{k\to\infty} \sum_{j\in\Gamma} \wei \omega_j^p |x_j^{(k)}|^p\leq \rho^p$ for all finite subsets $\Gamma \subset \Lambda$. This shows $\lnorm { x} {\wei \omega} p\leq \rho$.

In the case that $D_F$ is a ball centered at some reference solution  $x_0\neq 0$, we may replace the operator $F$ by the operator 
$x\mapsto F(x+x_0)$. This is equivalent to using the penalty term $\alpha \lnorm {x-x_0}{\wei r}1$ in \eqref{eq:Tik} with the original operator $F$, 
i.e.\ Tikhonov regularization with initial guess $x_0$. Without such a shift, \Cref{ass_domain} is violated.}

\subsection{Existence and uniqueness of minimizers}
We briefly address the question of existence and uniqueness of minimizers in \eqref{eq:Tik}. 
Existence follows by a standard argument of the direct method of the calculus of variations as often used in Tikhonov regularization, 
see, e.g., \cite[Thm.~3.22]{Scherzer_etal:09}).
\begin{prop} \label{prop:existence}
Suppose \Cref{ass_operator} holds true. Then for every $\gobs \in \Yspace$ and $\alpha>0$ there exists a solution to the minimization problem in \eqref{eq:Tik}.  
If $D_F= \lspace {\wei a} 2$ and $F$ is linear, then the minimizer is unique. 
\end{prop}
\begin{proof}
Let $(x^{(n)})_{n\in\mathbb{N}}$ be a minimizing sequence of the Tikhonov 
functional. Then $\lnorm{x^{(n)}}{\wei r}1$ is bounded. The compactness of the embedding $\lspace {\wei r}1\subset \lspace {\wei a}2$ (see \Cref{app:embed} in the appendix) implies the existence of a subsequence (w.l.o.g.\ again the full sequence) converging in \new{ $\lnorm {\cdot } {\wei a} 2$}
to some $x\in \lspace{\wei a}2$. Then $x\in D_F$ as $D_F$ is closed. The second inequality in Assumption \ref{ass_operator} implies 
\[\lim_{n\to\infty}\|g^{\mathrm{obs}}-F(x^{(n)})\|_{\Yspace}^2 = \|g^{\mathrm{obs}}-F(x)\|_{\Yspace}^2.\]
Moreover, for any finite subset $\Gamma\subset\Lambda$ we have 
\[ \sum_{j\in \Gamma } \wei r_j |x_j| = \lim_n  \sum_{j\in \Gamma } \wei r_j |x^{(n)}_j| 
\leq \liminf_n \lnorm{x^{(n)}}{\wei r}1,
\]
and hence $\lnorm{x}{\wei r}1\leq \liminf_n \lnorm{x^{(n)}}{\wei r}1$. This shows that $x$ minimizes the \linebreak
Tikhonov functional. 

%Existence follows from \cite[Thm.~3.22]{Scherzer_etal:09}). The statements $1.$,$2.$,$3.$ and $5.$ in \cite[Ass.~3.13]{Scherzer_etal:09}) are obviously satisfied  if we choose
%$U=\lspace {\wei r}1$ with $\tau_U$ the topology generated $\lnorm{\cdot}{\wei a}2$, $V=\Yspace$ with $\tau_V$ the norm topology, 
%and $\mathcal{R}(x)=\lnorm x {\wei r}1$. 
%To verify $4.$, suppose that $\lim_{n\to\infty}\lnorm{x^{(n)}-x} {\wei a}2=0$ and $\Gamma\subset\Lambda$ be a finite subset. Then 
%\[ \sum_{j\in \Gamma } \wei r_j |x_j| = \lim_n  \sum_{j\in \Gamma } \wei r_j |x^{(n)}_j| \leq \liminf_n \mathcal{R}(x^{(n)}).\]
%Hence $\mathcal{R}(x)\leq \liminf_n \mathcal{R}\left(x^{(n)}\right)$. The compactness of the embedding $\lspace {\wei r}1\subset \lspace {\wei a}2$ (see \Cref{app:embed}) yields $6.$, and $7.$ follows from the second inequality in Ass.~\ref{ass_operator} together with $4.$, $3.$ and the closedness of $D$. 
In the linear case the uniqueness follows from strict convexity.  
\end{proof}
\new{
Note that \Cref{prop:existence} also yields the existence of minimizers in \eqref{eq:Tik_besov} under Assumptions \ref{ass:wavelet} and \ref{ass_operator_Besov} and eqs.~\eqref{eqs:r_condition}. }

\new{ If $F=A\colon \lspace {\wei a} 2 \rightarrow \Yspace$ is linear and satisfies \Cref{ass_operator}, the usual argument (see, e.g., \cite[Lem.~ 2.1]{lorenz:08}) shows sparsity of the minimizers as follows: By the first order optimality condition there exists  $\xi\in\partial \lnorm{\cdot}{\wei r}1(\xh)$ such that $\xi$ belongs to the range of the adjoint $A^*$, that is $\xi\in \lspace {\wei a\inv}2$
and hence $\wei a_j\inv |\xi_j|\rightarrow 0$. Since $\wei a_j \wei r_j\inv \rightarrow 0$, we have  
$\wei a_j \leq \wei r_j$ for all but finitely many $j$. Hence, we obtain $|\xi_j|<r_j$, forcing $ x_j=0$ for all but finitely many $j$. \\ 
Note that for this argument to work, it is enough to require that $\wei a_j \wei r_j\inv$ is bounded from above. Also the existence of minimizers can be shown under this weaker assumption using  the weak$^\ast$-topology on $\lspace {\wei r} 1$ (see \cite[Prop.~2.2]{flemming:16}).} 

%Note that this suggests that for high dimensional problems a choice $r>0$ is required. \todo{e.g.\ for the Radon transform!}}
%\Cref{ass_operator} implies the existence and uniqueness of solutions to the minimization problem \eqref{eq:Tik} 
%by standard arguments (see, e.g.,  using the compactness of the embedding $\lspace {\wei r}1 \subset \lspace {\wei a}2$ shown in \Cref{app:embed}.

\section{Weak sequence spaces}\label{sec:weaksequencespace}
In this section we introduce spaces of sequences whose bounded sets will provide the source sets for the convergence analysis in the next chapters. We define a specific thresholding map and analyze its approximation properties. \\
Let us first introduce  a scale of spaces, part of which interpolates between the spaces $\lspace {\wei r}1$ and $\lspace {\wei a}2$ 
involved in our setting. 
For $t\in (0,2]$ we define weights  
\begin{align}\label{eq:omega}
(\wei \omega_t)_j:=(\wei a_j^{2t-2}\wei r_j^{2-t})^\frac{1}{t}. 
\end{align}
Note that $\wei \omega_1=\wei r$ and $\wei \omega_2=\wei a$. The next proposition captures interpolation inequalities we will need later.
\begin{prop}[Interpolation inequality] \label{prop:inter_scale}Let $u,v,t\in (0,2]$ and $\theta\in(0,1)$ with $\frac{1}{t}= \frac{1-\theta}{u} + \frac{\theta}{v}.$  
Then 
\[ \lnorm x{\wei \omega_t}t \leq \lnorm x {\wei \omega_{u}} {u} ^{1-\theta} \lnorm x {\wei \omega_{v}} {v} ^\theta \quad\text{for all } x\in \lspace {\wei \omega_u}u \cap \lspace {\wei \omega_v}v.\] 
\end{prop} 
\begin{proof}
We use Hölder's inequality with the conjugate exponents $\frac{u}{(1-\theta) t}$ and $\frac{v}{\theta t}$:
\begin{align*}
\lnorm x{\wei \omega_t}t^t 
& = \sum_{j\in\Lambda} \left( \wei a_j^{2u-2}  \wei r_j^{2-u} |x_j|^{u}\right)^\frac{(1-\theta) t}{u} \left( \wei a_j^{2v-2}  \wei r_j^{2-v} |x_j|^{v}\right)^\frac{\theta t}{v}\\ & \leq \lnorm x {\wei \omega_{u}} {u} ^{(1-\theta)t} \lnorm x {\wei \omega_{v}} {v} ^{\theta t}. 
\end{align*}
\end{proof}
\begin{rem}\label{rem:interpolation}
In the setting of \Cref{prop:inter_scale} real interpolation theory yields the stronger statement $\lspace {\wei \omega_t} t = (\lspace {\wei \omega_u}u , \lspace {\wei \omega_v}v )_{\theta,t}$ with equivalent quasi-norms (see, e.g., \cite[Theorem 2]{F:78}). The stated interpolation inequality is a consequence.
\end{rem} 
For $t\in (0,2)$ we define a weak version of the space $\lspace {\wei \omega_t}t$.
\begin{defi}[Source sets]\label{defi:kt} 
Let $t\in (0,2)$. We define \[k_t:= \{ x\in \mathbb{R}^\Lambda \colon \| x\|_{k_t}<\infty \}\]  with 
\[  \| x\|_{k_t}:= \sup_{\alpha>0} \alpha \left( \sum_{j\in \Lambda}  \wei a_j^{-2}\wei r_j^2 \mathds{1}_{ \{\wei a_j^{-2}\wei r_j  \alpha <|x_j| \} } \right)^\frac{1}{t}. \] 
\end{defi}
\begin{rem}\label{rem:weak}
The functions $ \| \cdot \|_{k_t}$ are quasi-norms. The quasi-Banach spaces $k_t$ are weighted Lorentz spaces. They appear as real interpolation spaces between weighted $L^p$ spaces. To be more precise  \cite[Theorem 2]{F:78} yields 
$ k_t= (\lspace {\wei \omega_u}u , \lspace {\wei \omega_v}v )_{\theta,\infty}$ with equivalence of quasi-norms for $u,v,t$ and $\theta$  as in \Cref{prop:inter_scale}.
\end{rem} 
%In \Cref{app:quasinorm} we show that $\varrho_t$ is a quasinorm on $k_t$. 
\begin{rem} \label{rem:markov}
\Cref{rem:interpolation} and \Cref{rem:weak}  predict an embedding \[ \lspace {\wei \omega_t} t =  (\lspace {\wei \omega_u}u , \lspace {\wei \omega_v}v )_{\theta,t} \subset  (\lspace {\wei \omega_u}u , \lspace {\wei \omega_v}v )_{\theta,\infty} = k_t.\]
Indeed the Markov-type inequality 
\[ \alpha^t \sum_{j\in \Lambda}\wei a_j^{-2} \wei r_j^2  \mathds{1}_{ \{ \wei a_j^{-2} \wei r_j \alpha < |x_j| \} }  \leq  \sum_{j\in \Lambda} \wei a_j^{2t-2} \wei r_j^{2-t} |x_j|^t = \lnorm x{\wei \omega_t} t^t \]   
proves $\| \cdot \|_{k_t} \leq \lnorm \cdot {\wei \omega_t} t $.
\end{rem} 
%\begin{prop}[Markov embedding] \label{prop:markov}
%Let $0<t<2$. Then $\varrho_t(\cdot)\leq  \lnorm \cdot{\wei \omega_t}t$. Hence, $\lspace {\wei \omega_t}t \subseteq k_t$. 
%\end{prop} 
%\begin{proof}
%The claim follows from  
%
%\end{proof}
For $\wei a_j=\wei r_j=1$ we obtain the weak $\ell_p$-spaces $k_t=\ell_{t,\infty}$ that appear in nonlinear approximation theory (see e.g. \cite{C:00}, \cite{C:09}).\\
We finish this section by defining a specific nonlinear thresholding procedure depending on $r$ and $a$ whose approximation theory is characterized by the spaces $k_t$. This characterization is the core for the proofs in the following chapters. 
The statement is \cite[Theorem 7.1]{C:00} for weighted sequence space. For sake of completeness we present an elementary proof based on a partition trick that is perceivable in the proof of \cite[Theorem 4.2]{C:00}.\\ 
Let $\alpha>0$. We consider the map 
\begin{align*}
 T_\alpha\colon \mathbb{R}^\Lambda \rightarrow \mathbb{R}^\Lambda \quad\text{by}\quad 
  T_\alpha(x)_j:= 
 \begin{cases} 
x_j  & \text{if } \wei a_j^{-2} \wei r_j \alpha < |x_j| \\
0 & \text{else } 
\end{cases} .
\end{align*} 
\new {Note that
\[ \alpha^2 \sum_{j\in \Lambda}  \wei a_j^{-2}\wei r_j^{2} \mathds{1}_{ \{ \wei a_j^{-2}\wei r_j \alpha < |x_j|\}} \leq \lnorm {T_\alpha(x)} {\wei a}2 ^2 \leq \lnorm {x} {\wei a}2 ^2. \]
If $\wei a_j \wei r_j\inv$ is bounded above, then $\wei a_j^{-2}\wei r_j^{2}$ is bounded away from zero. Hence, in this case we see that the set of $j\in \Lambda$ with $\wei a_j^{-2}\wei r_j \alpha < |x_j|$ is finite, i.e.\ $T_\alpha(x)$ has only finitely many nonvanishing coefficients whenever $x\in \lspace {\wei a}2$.}
\begin{lem}[Approximation rates for $T_\alpha$] \label{threshold_rates} 
Let $0<t<p\leq 2$ and $x\in \mathbb{R}^\Lambda$. Then $x\in k_t$ if and only if 
$\gamma(x):= \sup_{\alpha>0} \alpha^\frac{t-p}{p}\lnorm {x-T_\alpha(x)} {\wei \omega_p} p < \infty $.\\
More precisely we show bounds \[ \gamma(x)\leq 2\left(2^{p-t}-1\right)^{-\frac{1}{p}}  \| x\|_{k_t}^\frac{t}{p} \quad\text{and}\quad  \| x\|_{k_t}\leq 2^\frac{p}{t}(2^t-1)^{-\frac{1}{t}}\gamma(x)^\frac{p}{t}. \]
\end{lem} 
\begin{proof} 
We use a partitioning to estimate 
\begin{align*}
\lnorm {x-T_\alpha(x)} {\wei \omega_p} p ^p & = \sum_{j\in \Lambda}  \wei a_j^{2p-2}\wei r_j^{2-p} |x_j|^p \mathds{1}_{ \{ |x_j|\leq \wei a_j^{-2}\wei r_j \alpha\}} \\
&= \sum_{k=0}^\infty \sum_{j\in \Lambda} \wei a_j^{2p-2}\wei r_j^{2-p} |x_j|^p  \mathds{1}_{ \{ \wei a_j^{-2}\wei r_j 2^{-(k+1)}\alpha  <|x_j|\leq \wei a_j^{-2}\wei r_j 2^{-k}\alpha\}}\\ 
&\leq \alpha^p \sum_{k=0}^\infty 2^{-pk} \sum_{j\in \Lambda}  \wei a_j^{-2}\wei r_j^{2}  \mathds{1}_{ \{ \wei a_j^{-2}\wei r_j 2^{-(k+1)}\alpha <|x_j|\}} \\
&\leq  \alpha^{p-t}  \| x\|_{k_t}^t  2^t  \sum_{k=0}^\infty (2^{t-p})^k\\
&= \alpha^{p-t} 2^p \left(2^{p-t}-1\right)\inv   \| x\|_{k_t}^t.
\end{align*}
A similar estimation yields the second inequality: 
\begin{align*}
\sum_{j\in \Lambda}  \wei a_j^{-2}\wei r_j^{2}  \mathds{1}_{ \{ \wei a_j^{-2} \wei r_j\alpha < |x_j|\}} &=  \sum_{k=0}^\infty \sum_{j\in \Lambda}  \wei a_j^{-2}\wei r_j^{2}  \mathds{1}_{ \{ \wei a_j^{-2}\wei r_j 2^k\alpha  < |x_j|\leq \wei a_j^{-2} \wei r_j 2^{k+1}\alpha   \}} \\
& \leq \alpha^{-p} \sum_{k=0}^\infty 2^{-kp}  \sum_{j\in \Lambda}  \wei a_j^{2p-2}\wei r_j^{2-p} |x_j|^p  \mathds{1}_{\{|x_j|\leq \wei a_j^{-2} \wei r_j 2^{k+1}\alpha   \}}\\&= \alpha^{-p} \sum_{k=0}^\infty 2^{-kp} \lnorm {x-T_{2^{k+1}\alpha }(x)}  {\wei \omega_p} p ^p \\
& \leq \alpha^{-t} \gamma(x)^p 2^{p-t}  \sum_{k=0}^\infty (2^{-t})^k \\ 
&= \alpha^{-t}  \gamma(x)^p 2^p \left(2^t-1\right)\inv.
\end{align*}
\end{proof}

\begin{cor}\label{cor:embedweak}
Assume $\wei a_j\wei r_j\inv$ is bounded from above.
Let $0<t<p\leq 2$. Then $k_t\subset \lspace {\wei \omega_p}p$. More precisely, there is a constant $M>0$ depending on $t,p$ and $\sup_{j\in\Lambda}\wei a_j\wei r_j\inv$ such that $ \lnorm \cdot {\wei \omega_p} p\leq M   \| \cdot\|_{k_t}$. 
\end{cor} 
\begin{proof}
Let $x\in k_t$. The assumption implies the existence of a constant $c>0$ with $c\leq \wei a_j^{-2} \wei r_j^2$ for all $j\in \Lambda.$ Let $\alpha>0$. Then
\[ c \sum_{j\in\Lambda}  \mathds{1}_{ \{\wei a_j^{-2}\wei r_j  \alpha <|x_j| \} } \leq  \sum_{j\in\Lambda}  \wei a_j^{-2}\wei r_j^2 \mathds{1}_{ \{\wei a_j^{-2}\wei r_j  \alpha < |x_j| \} }\leq  \| x\|_{k_t}^t \alpha^{-t}.\] Inserting $\overline{\alpha}:= 2  \| x\|_{k_t} c^{-\frac{1}{t}}$ implies $\wei a_j^{-2}\wei r_j \overline{\alpha}\geq |x_j|$ for all $j\in \Lambda.$ Hence, $T_{\overline{\alpha}}(x)=0$. With ${C=2\left(2^{p-t}-1\right)^{-\frac{1}{p}}}$ 
 \Cref{threshold_rates} yields
 \[\lnorm x {\wei \omega_p} p= \lnorm {x-T_{\overline{\alpha} }(x)} {\wei \omega_p} p\leq C \|x\|_{k_t}^\frac{t}{p} \overline{\alpha}^\frac{p-t}{p}=2^\frac{p-t}{p}C c^\frac{t-p}{tp}  \|x\|_{k_t}.  \]    
\end{proof}

\new{
\begin{rem}[Connection to best $N$- term approximation]
For better understanding of the source sets we sketch another characterization of $k_t$. 
For $z\in \mathbb{R}^\Lambda$ we set $S(x):= \sum_{j\in \Lambda} \wei a_j^{-2} \wei r_j^2 \mathds{1}_{ \{z_j \neq 0\} }.$ Note that for $\wei a_j=\wei r_j=1$ we simply have $S(x)= \# \mathrm{supp}(x)$. Then for $N>0$ one defines the best approximation error by 
\[ \sigma_N(x):= \inf\left\{ \lnorm{x-z}{\wei a}2 \colon S(z)\leq N\right\}. \] 
Using arguments similar to those in the proof of $\Cref{lem:threshold_bounds}$ one can show that for $t \in (0,2)$ 
we have $x\in k_t$ if and only if the error scales like $\sigma_N(x)=\mathcal{O}(N^{\frac{1}{2}-\frac{1}{t}})$. 
\end{rem}}
\section{Convergence Rates via Variational Source Conditions}\label{sec:non_oversmoothing}
We prove rates of convergence for the regularization scheme \eqref{eq:Tik} based on variational source conditions. The latter are nessecary and often sufficient conditions for rates of convergence for Tikhonov regularization and other regularization methods (\cite{Scherzer_etal:09,flemming:12b,HW:17}). For $\ell^1$-norms these conditions are typically of the form
\begin{align}\label{eq:vsc}
\beta \lnorm {x^\dagger - x}{\wei r}1 + \lnorm {x^\dagger}{\wei r}1 - \lnorm {x}{\wei r}1 \leq \psi\left(\|F(x)- F(x^\dagger)\|_\Yspace^2 \right) \quad\text{for all } x \in D_F\cap \lspace {\wei r}1
\end{align} 
with $\beta\in [0,1]$ and $\psi\colon [0,\infty)\rightarrow [0,\infty)$ a concave, stricly increasing function with $\psi(0)=0$.  
%In \cite{burger2013convergence} the authors verify \eqref{eq:vsc} with $\beta=1$ for a linear operator under the assumption that the unit vectors lie in the range of the adjoint operator. \todo{Dopplung mit Einleitung vermeiden?} The latter condition can be interpreted as a smoothness condition on the basis functions. A relaxation of this assumption is studied in \cite{FH:15} and leads to \eqref{eq:vsc} for $\beta\in (0,1]$. In \cite{flemming2018injectivity} the authors show that these relaxed assumptions are satisfied under mild assumptions and conclude that  \eqref{eq:vsc} is always statisfied for an appropriate index function depending on decay of the sequence $x^\dagger.$  
The common starting point of verifications of \eqref{eq:vsc} in the references \cite{burger2013convergence,FH:15,flemming2018injectivity,HM:19}, which  
have already been discussed in the introduction,   
is a splitting of the left hand side in \eqref{eq:vsc} into two summands according 
to a partition of the index set into low level and high level indices. 
The key difference to our verification in \cite{HM:19} is that this partition will be chosen adaptively to $x^\dagger$ below. 
This possibility is already mentioned, but not further exploited in \cite[Remark 2.4]{Flemming2016} and 
\cite[Chapter 5]{flemming2018injectivity}.

\subsection{variational source conditions}
We start with a Bernstein-type inequality. 
\begin{lem}[Bernstein inequality]\label{lemma:bernstein}
Let $t\in  (0,2)$, $x^\dagger\in k_t$ and $\alpha>0$. We consider \[\Lambda_\alpha:=\{j\in\Lambda \colon \wei a_j^{-2}\wei r_j  \alpha < |x_j^\dagger| \}\] and the coordinate projection $P_\alpha \colon \mathbb{R}^\Lambda \rightarrow  \mathbb{R}^\Lambda$ onto $\Lambda_\alpha$ given by $(P_\alpha x)_j:= x_j$ if $j\in \Lambda_\alpha$ and ${(P_\alpha x)_j:= 0}$ else. Then  
\[ \lnorm {P_\alpha x} {\wei r}1 \leq \| x^\dagger\|_{k_t}^\frac{t}{2} \alpha^{-\frac{t}{2}} \lnorm x {\wei a} 2 \quad\text{for all} \quad x\in \lspace {\wei a}2.\]
\end{lem}
\begin{proof}
Using the Cauchy–Schwarz inequality we obtain
\begin{align*}
 \lnorm {P_\alpha x} {\wei r}1& = \sum_{j\in\Lambda}  \left(  \wei a_j\inv  \wei r_j \mathds{1}_{ \{\wei a_j^{-2}\wei r_j  \alpha < |x_j^\dagger| \} } \right) \left(\wei a_j |x_j|\right)  \\
 &\leq \left(   \sum_{j\in\Lambda} \wei a_j^{-2}\wei r_j^2  \mathds{1}_{ \{\wei a_j^{-2}\wei r_j  \alpha < |x_j^\dagger |\} } \right)^\frac{1}{2} \left(  \sum_{j\in\Lambda}\wei a_j^2 |x_j|^2 \right)^\frac{1}{2} \\
&\leq \| x^\dagger\|_{k_t}^\frac{t}{2} \alpha^{-\frac{t}{2}}  \lnorm x {\wei a}2.
\end{align*}
\end{proof}
The following lemma characterizes variational source conditions \eqref{eq:vsc} for the embedding operator $\ell^1_r\hookrightarrow \ell^2_a$ (if $\wei a_j \wei r_j\inv \rightarrow 0$)
and power-type functions $\psi$ 
with $\beta=1$ and $\beta=0$ in terms of the weak sequence spaces $k_t$ in \Cref{defi:kt}:
\begin{lem}[variational source condition for embedding operator]\label{lemma:vsc_id} 
\item Assume $x^\dagger\in\lspace {\wei r}1$ and $t\in (0,1)$. 
The following statements are equivalent: 
\begin{enumerate}[(i)]
\item $x^\dagger\in k_t.$  
\item There exist a constant $K>0$  such that 
\begin{align}\label{eq:VSC}
 \lnorm {x^\dagger-x}{\wei r}1 + \lnorm {x^\dagger}{\wei r}1 -\lnorm x{\wei r}1\leq K \lnorm{x^\dagger-x}{\wei a}2 ^\frac{2-2t}{2-t}
\end{align}
for all  $x\in  \lspace {\wei r}1.$
\item There exist a constant $K>0$ such that 
\[ \lnorm {x^\dagger}{\wei r}1 -\lnorm x{\wei r}1\leq K \lnorm{x^\dagger-x}{\wei a}2 ^\frac{2-2t}{2-t}\]
for all $x\in \lspace {\wei r}1$ with  $|x_j| \leq |x^\dagger_j|$ for all $ j\in \Lambda.$
\end{enumerate}
More precisely, $(i)$ implies $(ii)$ with $K=(2+4(2^{1-t}-1)^{-1}) \|x^\dagger\|_{k_t}^\frac{t}{2-t}$ and $(iii)$ yields the bound $\|x^\dagger\|_{k_t}\leq K^\frac{2-t}{t}.$
\end{lem}
\begin{proof}
First we assume $(i)$. For $\alpha>0$ we consider $P_\alpha$ as defined in \Cref{lemma:bernstein}. Let $x\in D\cap \lspace {\wei r}1$. By splitting all three norm term in the left hand side of \eqref{eq:VSC} by $\lnorm {\cdot}{\wei r}1=\lnorm{P_\alpha \cdot} {\wei r}1+\lnorm {(I-P_\alpha )\cdot}{\wei r}1$ and using the triangle equality for the $(I-P_\alpha)$ terms and the reverse triangle inequality for the $P_\alpha$ terms (see \cite[Lemma 5.1] {burger2013convergence}) we obtain 
\begin{align}\label{eq:ansatz_vsc}
\lnorm {x^\dagger-x}{\wei r}1 + \lnorm {x^\dagger}{\wei r}1   -\lnorm x{\wei r}1 \leq 2\lnorm {P_\alpha(x^\dagger-x)}{\wei r}1 + 2 \lnorm {(I-P_\alpha) x^\dagger}{\wei r}1.
\end{align}
We use \Cref{lemma:bernstein} to handle the first summand
\[ \lnorm {P_\alpha(x^\dagger-x)}{\wei r}1 \leq  \| x^\dagger\|_{k_t}^\frac{t}{2} \alpha^{-\frac{t}{2}} \lnorm {x^\dagger-x} a 2.\] 
Note that $P_\alpha x^\dagger = T_\alpha( x^\dagger).$ Hence,  
\Cref{threshold_rates} yields %a constant $c$ that depends only on $t$ such that 
\[  \lnorm {(I-P_\alpha) x^\dagger}{\wei r}1 = \lnorm {x^\dagger-T_\alpha( x^\dagger)}{\wei r}1
\leq 2(2^{1-t}-1)^{-1}\| x^\dagger\|_{k_t}^t \alpha^{1-t}.\]
Inserting the last two inequalities into \eqref{eq:ansatz_vsc} and choosing \[ \alpha=  \lnorm {x^\dagger-x} a 2^\frac{2}{2-t} \| x^\dagger\|_{k_t}^{-\frac{t}{2-t}}\] we get $(ii)$. 
%we obtain 
%\begin{align} \label{eq:vsc_replacement}
% \lnorm {x^\dagger-x}{\wei r}1 + \lnorm {x^\dagger}{\wei r}1   -\lnorm x{\wei r}1 
% \leq  2\| x^\dagger\|_{k_t}^\frac{t}{2} \alpha^{-\frac{t}{2}} \lnorm {x^\dagger-x} a 2 + 4(2^{1-t}-1)^{-1}\| x^\dagger\|_{k_t}^t \alpha^{1-t} 
%\end{align}
%for all $\alpha>0$. 
%Choosing $\alpha=  \lnorm {x^\dagger-x} a 2^\frac{2}{2-t} \| x^\dagger\|_{k_t}^{-\frac{t}{2-t}}$ we get $(ii)$.
% with $K=(c+1)\varrho_t(x^\dagger)^{\frac{t}{2-t}}$.
\\ 
Obviously $(ii)$ implies $(iii)$ as $\lnorm {x^\dagger-x}{\wei r}1\geq 0.$\\ 
It remains to show that $(iii)$ implies $(i)$. Let $\alpha>0$. We define 
\[ x_j:=  \begin{cases} 
x_j^\dagger  & \text{if } |x^\dagger_j| \leq  \wei a_j^{-2} \wei r_j \alpha   \\
x_j^\dagger- \wei a_j^{-2} \wei r_j \alpha  & \text{if } x^\dagger_j > \wei a_j^{-2} \wei r_j \alpha  \\
x_j^\dagger+ \wei a_j^{-2} \wei r_j \alpha  & \text{if } x^\dagger_j < - \wei a_j^{-2} \wei r_j \alpha  
\end{cases} .\] 
Then $|x_j|\leq |x^\dagger_j|$ for all $j\in\Lambda$. Hence, $x\in \lspace {\wei r}1$. We estimate 
\begin{align*}
\alpha \sum_{j\in \Lambda}  \wei a_j^{-2}\wei r_j^2 \mathds{1}_{ \{\wei a_j^{-2}\wei r_j  \alpha  <|x^\dagger_j| \} } & = \lnorm {x^\dagger}{\wei r}1 -\lnorm x{\wei r}1 
    \leq K \lnorm{x^\dagger-x}{\wei a}2 ^\frac{2-2t}{2-t} \\ 
   &= K \left( \sum_{j\in \Lambda} \wei a_j^2 (\wei a_j^{-2} \wei r_j \alpha)^2 \mathds{1}_{ \{\wei a_j^{-2}\wei r_j  \alpha  < |x^\dagger_j| \} } \right)^\frac{1-t}{2-t} \\ 
   &= K \alpha^\frac{2-2t}{2-t} \left( \sum_{j\in \Lambda}  \wei a_j^{-2}\wei r_j^2 \mathds{1}_{ \{\wei a_j^{-2}\wei r_j  \alpha  < |x^\dagger_j| \} } \right)^\frac{1-t}{2-t}.
\end{align*}
Rearranging terms in this inequality yields  
\[ \sum_{j\in \Lambda}  \wei a_j^{-2}\wei r_j^2 \mathds{1}_{ \{\wei a_j^{-2}\wei r_j  \alpha < |x^\dagger_j| \} }\leq K^{2-t} \alpha^{-t}.\]
Hence, $ \|x^\dagger\|_{k_t}\leq K^\frac{2-t}{t}.$  
\end{proof}

\begin{thm}[variational source condition]\label{vsc_nonlinear}
Suppose \Cref{ass_operator} holds true and let $t\in(0,1)$, $\varrho>0$ and $x^\dagger\in D$. If $\|x\|_{k_t}\leq \varrho$ 
then the variational source condition  
\begin{align}\label{eq:vsc_nonlinear}\nonumber
&\lnorm {x^\dagger-x}{\wei r}1 + \lnorm {x^\dagger}{\wei r}1 -\lnorm x{\wei r}1\leq C_\mathrm{vsc}  \| F(x^\dagger)-F(x) \|_\Yspace^\frac{2-2t}{2-t}\\
& \mbox{for all } x\in  D_F\cap \ell^r_1
\end{align}  
holds true with $C_\mathrm{vsc} = (2+4(2^{1-t}-1)^{-1}) L^{\frac{2-2t}{2-t}}\varrho^\frac{t}{2-t}$.\\
If in addition \Cref{ass_domain} holds true, then \eqref{eq:vsc_nonlinear} implies $\|x\|_{k_t}\leq L^{\frac{2-2t}{t}} C_\mathrm{vsc}^{\frac{2-t}{t}}$.
\end{thm}
\begin{proof}
\Cref{cor:embedweak} implies $x\in D\cap \lspace {\wei r}1$. The first claim follows from the first inequality in \Cref{ass_operator} together with \Cref{lemma:vsc_id}. The second inequality in \Cref{ass_operator} together with \Cref{ass_domain} imply statement (iii) in \Cref{lemma:vsc_id} with $K= L^\frac{2-2t}{2-t} C_\mathrm{vsc}.$ Therefore, \Cref{lemma:vsc_id} yields the second claim.
\end{proof}

\subsection{Rates of Convergence}\label{sec:rates}
In this section we formulate and discuss bounds on the reconstruction error which follow from the variational source condition
\eqref{eq:vsc_nonlinear} by general variational regularization theory (see, e.g., \cite[Prop. 4.2, Thm. 4.3]{HM:19} or  
\cite[Prop.13., Prop.14.] {flemming2018injectivity}).
%
%For convenience of the reader we included a proof in  \cref{proof:rates} that is based only on \eqref{eq:vsc_replacement} and \Cref{ass_operator}.   

%\begin{prop}\label{prop:exact_rates}
%Let $x^\dagger\in D\cap \lspace {\wei r}1.$ Assume that the variational source condition \eqref{eq:vsc_nonlinear} holds true. Then the minimizers $x_\alpha$ of \eqref{eq:Tik} for exact data $\gobs=F(x^\dagger)$ satisfy the error bounds 
%\begin{align*}
%\|F(x^\dagger)-F(x_\alpha)\|_\Yspace & \leq (2C_{vsc})^\frac{2-t}{2} \varrho^\frac{t}{2} \alpha^\frac{2-t}{2}\quad \text{and } \\ 
%\lnorm{x^\dagger-x_\alpha}{\wei r}1 &  \leq 2^{1-t} C_{vsc}^{2-t} \varrho^t\alpha^{1-t} \quad\text{for all } \alpha>0
%\end{align*}
%\end{prop}  
%\begin{proof}
%As $x_\alpha$ is a minimizer of $\eqref{eq:Tik}$ and by \eqref{eq:vsc_nonlinear}  we have 
%\[ \frac{1}{2} \|F(x^\dagger)-F(x_\alpha)\|_\Yspace^2 \leq \alpha \left( \lnorm {x^\dagger}{\wei r}1 - \lnorm {x_\alpha}{\wei r}1  \right)\leq \alpha C_{vsc} \varrho^\frac{t}{2-t} \| F(x^\dagger)-F(x_\alpha) \|_\Yspace^\frac{2-2t}{2-t} .\] 
%Solving this inequality yields the first bound. Moreover, we obtain $\lnorm {x^\dagger}{\wei r}1 - \lnorm {x_\alpha}{\wei r}1 \geq 0$. Therefore, the second bound follows by inserting the first into \eqref{eq:vsc_nonlinear}. 
%\end{proof}
\begin{thm}[Convergence rates]\label{thm:rates}
Suppose \Cref{ass_operator} holds true. Let $t\in(0,1)$, $\varrho>0$ and $x^\dagger\in D_F$ with $\|x^\dagger\|_{k_t}\leq \varrho.$ Let $\delta\geq 0$ and $\gobs\in \Yspace$ satisfy $\|\gobs-F(x^\dagger)\|_\Yspace\leq\delta$.
\begin{enumerate}
\item (error splitting) 
Every minimizer $\xh$ of \eqref{eq:Tik} satisfies 
\begin{align}\label{eq:rate_in_r1}
\lnorm{x^\dagger-\xh}{\wei r}1 &  \leq C_e \left( \delta^2 \alpha\inv + \varrho^t \alpha^{1-t} \right) \quad\text{and} \\ \label{eq:rate_in_a2}
\lnorm {x^\dagger-\xh } a 2 
&\leq C_e \left(\delta +  \varrho^\frac{t}{2} \alpha^\frac{2-t}{2}\right).
\end{align}  
for all $\alpha>0$ with a constant $C_{e}$ depending only on $t$ and $L$.
\item  (rates with a-priori choice of $\alpha$)
If $\delta>0$ and $\alpha$ is chosen such that \[ c_1 \varrho^\frac{t}{t-2} \delta^\frac{2}{2-t} \leq  \alpha \leq c_2 \varrho^\frac{t}{t-2} \delta^\frac{2}{2-t} \quad\text{for } 0<c_1<c_2,\]  then every minimizer $\hat{x}_\alpha$ of \eqref{eq:Tik} satisfies 
\begin{align}\label{eq:rate_appri}
 \lnorm {x^\dagger-\xh} {\wei r}1  & \leq C_p \varrho^\frac{t}{2-t}\delta^\frac{2-2t}{2-t} \quad\text{and} \\
  \lnorm {x^\dagger-\xh} {\wei a}2 &\leq C_p \delta.
\end{align}  
with a constant $C_{p}$ depending only on $c_1, c_2, t$ and $L$.
\item  (rates with discrepancy principle) 
Let $1\leq\tau_1\leq \tau_2$. If $\hat{x}_\alpha$ is a minimizer of \eqref{eq:Tik} with $\tau_1 \delta\leq \|F(\hat{x}_\alpha)-\gobs \|_\Yspace \leq \tau_2 \delta$, then 
\begin{align}\label{eq:rate_discr}
\lnorm {x^\dagger-\xh} {\wei r}1 & \leq C_d \varrho^\frac{t}{2-t}\delta^\frac{2-2t}{2-t}   \quad\text{and} \\
\lnorm {x^\dagger-\xh} {\wei a}2 & \leq C_d  \delta. 
\end{align}
Here $C_d>0$ denotes a constant depending only on $\tau_2$, $t$ and $L$. 
\end{enumerate}
\end{thm}
%\begin{proof}
%By \Cref{vsc_nonlinear} the variational source condition \eqref{eq:vsc_nonlinear} is satisfied. Hence, \textit{1.} follows from \Cref{prop:exact_rates}. For proofs of \textit{2.} and \textit{3.} we refer to \cite[Prop.13., Prop.14.]{flemming2018injectivity}.
%\end{proof}

We discuss our results in the following series of remarks: 

\begin{rem}
The proof of \Cref{thm:rates} makes no use of the second inequality in \Cref{ass_operator}. 
\end{rem}
\begin{rem}[error bounds in intermediate norms]  \label{rem:error_bounds_intermediate}
Invoking the interpolation inequalities given in \Cref{prop:inter_scale} allows to combine the bounds in the norms $\lnorm \cdot {\wei r}1$ and $\lnorm \cdot {\wei a}2$ to bounds in $\lnorm \cdot {\wei \omega_p}p$ for $p\in (t,1]$. \\ 
In the setting of \Cref{thm:rates}$(2.)$ or $(3.)$ we obtain 
\begin{align} \label{eq:error_in_e}
 \lnorm {x^\dagger - \xh} {\wei \omega_p} p \leq C \varrho^{\frac{t}{p}\frac{2-p}{2-t}} \delta^{\frac{2}{p}\frac{p-t}{2-t}}
\end{align}
with $C=C_p$ or $C=C_d$ respectively. 
\end{rem} 
\begin{rem} [Limit $t\rightarrow 1$]\label{rem:tequal1}
Let us consider the limiting case $t=1$ by assuming only $x^\dagger\in \lspace {\wei r}1\cap D_F$. Then it is well known, that the parameter choice $\alpha\sim \delta^2$ as well the discrepancy principle as in \Cref{thm:rates}.$3.$ lead to  bounds $\lnorm {x^\dagger-\xh} {\wei r}1 \leq C \lnorm {x^\dagger} {\wei r}1 $ and $\|F(x^\dagger)-F(\xh)\|_\Yspace\leq C \delta$. As above, \Cref{ass_operator} allows to transfer to a bound $\lnorm {x^\dagger-\xh} {\wei a}2\leq \tilde{C} \delta.$ Interpolating as in the last remark yields \[ \lnorm {x^\dagger-\xh} {\wei \omega_p}p \leq \tilde{C} \lnorm {x^\dagger} {\wei r}1^\frac{2-p}{p} \delta^\frac{2p-2}{p}. \] 
%Later this will serve as a replacement for \eqref{eq:error_in_e} in the case $t=1$.  
\end{rem} 

\begin{rem}[Limit $t\to 0$]
Note that in the limit $t\to 0$ the convergence rates get arbitrarily close to the linear convergence rate 
$\mathcal{O}(\delta)$, i.e., in contrast to standard quadratic Tikhonov regularization in Hilbert spaces no saturation 
effect occurs. This is also the reason why we always obtain optimal rates with the discrepancy principle even for smooth solutions $x^\dagger$.  

As already mentioned in the introduction, the formal limiting rate for $t\to 0$, i.e.\ a linear convergence rate in $\delta$ occurs if 
and only if $x^\dagger$ is sparse as shown by different methods in \cite{GHS:11}.
\end{rem}
We finish this subsection by showing that the convergence rates \eqref{eq:rate_appri},  \eqref{eq:rate_discr}, 
and \eqref{eq:error_in_e} are optimal in a minimax sense. 
\begin{prop}[Optimality]\label{prop:optimal}
Suppose that \Cref{ass_operator} holds true. Assume furthermore that there are $c_0>0$, $q\in (0,1)$ such that for every $\eta\in (0,c_0)$ there is $j\in \Lambda$ satisfying ${q\eta \leq \wei a_j \wei r_j\inv \leq \eta}$. Let $p\in (0,2]$, $t\in (0,p)$ and $\rho>0$. Suppose  $D$ contains all $x\in k_t$ with $\|x\|_{k_t}\leq \varrho.$ Consider an arbitrary reconstruction 
method described by a mapping $R:\Yspace\to \ell^1_r$ approximating the inverse of $F$. Then the worst case error under 
the a-priori information $\|x^\dagger\|_{k_t}\leq \varrho$ is bounded below by 
\begin{align}\label{eq:lower_bound}\nonumber
\sup \left\{ \lnorm{R\left(g^{\mathrm{obs}}\right)-x^\dagger}{\wei \omega_p} p: \left\|x^\dagger\right\|_{k_t}\leq \rho, 
\|F(x^\dagger)-g^{\mathrm{obs}}\|_{\Yspace}\leq \delta \right\}\\ 
\geq c \varrho^{\frac{t}{p}\frac{2-p}{2-t}} \delta^{\frac{2}{p}\frac{p-t}{2-t}}.
\end{align}
for all $\delta\leq \frac{1}{2}L\varrho c_0^\frac{2-t}{t}$
with $c=q^\frac{2p-2t}{pt} (2L\inv)^{\frac{2}{p}\frac{p-t}{2-t}}$. 
\end{prop}

\begin{proof}
It is a well-known fact that the left hand side in \eqref{eq:lower_bound} is bounded from below by $\frac{1}{2}\Omega(2\delta ,\varrho)$ with the modulus of continuity 
\begin{align*}
&\Omega(\delta,\varrho):= \\& \sup\left\{\lnorm{x^{(1)}-x^{(2)}}{\wei \omega_p} p: \left\|x^{(1)}\right\|_{k_t}, 
\left\|x^{(2)}\right\|_{k_t}\leq\rho, 
\left\| F\left(x^{(1)}\right)-F\left(x^{(2)}\right)\right\|_{\Yspace}\leq \delta \right\} 
\end{align*}
(see {\cite[Rem. 3.12]{EHN:96}}, {\cite[Lemma 2.8]{WSH:20}}). By Assumption \ref{ass_operator} we have 
\[
\Omega(\delta,\rho)\geq \sup\{\|x\|_{\wei \omega_p, p}: \|x\|_{k_t}\leq \rho, \|x\|_{\wei a,2}\leq 2L\inv \delta\}.
\]
By assumption there exists $j_0\in \Lambda$ such that 
\[ q\left(2L\inv \delta\varrho\inv \right)^\frac{t}{2-t}\leq \wei a_{j_0} \wei r_{j_0}\inv \leq  \left(2L\inv \delta\varrho\inv \right)^\frac{t}{2-t}. \] 
Choosing $x_{j_0}=\varrho a_{j_0}^\frac{2-2t}{t}r_{j_0}^\frac{t-2}{t}$ and  $x_j=0$ if $j\neq j_0$ we obtain 
$ \|x\|_{k_t}=\varrho $  and $ \lnorm x {\wei a}2 \leq 2 L\inv \delta$ and estimate 
\[ \lnorm x {\wei \omega_p} p = \varrho \left(\wei a_{j_0} \wei r_{j_0}\inv\right)^\frac{2p-2t}{pt}\geq q^\frac{2p-2t}{pt} (2L\inv)^{\frac{2}{p}\frac{p-t}{2-t}}  \varrho^{\frac{t}{p}\frac{2-p}{2-t}} \delta^{\frac{2}{p}\frac{p-t}{2-t}}.   \] 
\end{proof}
Note that for $\Lambda=\mathbb{N}$ the additional assumption in \Cref{prop:optimal} is satisfied if 
$\wei a_j \wei r_j\inv\sim \tilde{q}^j$ for $\tilde{q}\in (0,1)$ or if 
$\wei a_j \wei r_j\inv\sim j^{-\kappa}$ for $\kappa>0$, but violated if 
$\wei a_j \wei r_j\inv\sim \exp(-j^2)$.

\subsection{Converse Result} 
As a main result, we now prove that the condition $x^\dagger\in k_t$ is necessary and sufficient for the 
H\"older type approximation rate $\mathcal{O}(\alpha^{1-t})$:
\begin{thm}[converse result for exact data] \label{thm:converse}
Suppose \Cref{ass_operator} and \ref{ass_domain} hold true. Let $x^\dagger\in D_F\cap \lspace {\wei r}1$, $t\in (0,1)$, 
and $(x_\alpha)_{\alpha>0}$ the minimizers of \eqref{eq:Tik} for exact data $\gobs = F(x^\dagger).$ Then the following statements are equivalent: 
\begin{enumerate}[(i)]
\item $x^\dagger \in k_t.$
\item There exists a constant $C_2>0$ such that $\lnorm {x^\dagger-x_\alpha}{\wei r}1\leq C_2 \alpha^{1-t}$ for all $\alpha>0$. 
\item There exists a constant $C_3>0$ such that $\|F(x^\dagger)-F(x_\alpha)\|_\Yspace\leq C_3 \alpha^\frac{2-t}{2}$ for all $\alpha>0.$
\end{enumerate}   
More precisely, we can choose $C_2:= c \|x^\dagger\|_{k_t}^t$, $C_3:= \sqrt{2C_2}$ and bound \linebreak $\|x^\dagger\|_{k_t}\leq c C_3^\frac{2}{t}$ with a constant $c>0$ that depends on $L$ and $t$ only. 
\end{thm}
\begin{proof}
\item[$(i)\Rightarrow (ii)$:] By \Cref{thm:rates}(1.) for $\delta=0.$
\item[$(ii)\Rightarrow (iii)$:] As $x_\alpha$ is a minimizer of \eqref{eq:Tik} we have 
\[ \frac{1}{2} \|F(x^\dagger)-F(x_\alpha)\|_\Yspace^2 \leq \alpha\left( \lnorm {x^\dagger} {\wei r}1 - \lnorm {x_\alpha}{\wei r}1\right) \leq \alpha \lnorm {x^\dagger-x_\alpha}{\wei r}1\leq C_2 \alpha^{2-t}. \] 
Multiplying by $2$ and taking square roots on both sides yields $(iii)$. 
 \item[$(iii)\Rightarrow (i)$:] The strategy is to prove that $\|F(x^\dagger)-F(x_\alpha)\|_\Yspace$ is an upper bound on $\lnorm{x^\dagger-T_\alpha(x^\dagger)}{\wei a}2$ up to a constant and a linear change of $\alpha$ and then proceed using \Cref{threshold_rates}. \\
As an intermediate step we first consider  
\begin{align}\label{eq:approximation_for_converse}
 z_\alpha \in \argmin_{z\in \lspace {\wei r}1 } \left( \frac{1}{2} \lnorm {x^\dagger -z}{\wei a}2 ^2 +\alpha \lnorm{z}{\wei r}1 \right).
\end{align} 
The minimizer can be calculated in each coordinate separately by 
\begin{align*}
(z_\alpha)_j & =\argmin_{z\in \mathbb{R}} \left( \frac{1}{2} \wei a_j^2 |x^\dagger_j-z|^2 + \alpha \wei r_j |z| \right) \\ & = \argmin_{z\in \mathbb{R}} \left( \frac{1}{2} |x^\dagger_j-z|^2 + \alpha \wei a_j^{-2} \wei r_j |z| \right).
\end{align*}
Hence,  
\begin{align*}
 (z_\alpha)_j=  \begin{cases} 
x_j^\dagger-\wei a_j^{-2} \wei r_j \alpha  & \text{if } x^\dagger_j > \wei a_j^{-2} \wei r_j \alpha   \\
x_j^\dagger+ \wei a_j^{-2} \wei r_j \alpha  & \text{if } x^\dagger_j < -\wei a_j^{-2} \wei r_j \alpha  \\
0  & \text{if } |x^\dagger_j| \leq - \wei a_j^{-2} \wei r_j \alpha  
\end{cases}. 
\end{align*}
Comparing  $z_\alpha$ with $T_\alpha(x^\dagger)$ yields 
$ |x^\dagger-T_\alpha(x^\dagger)_j|\leq |x^\dagger_j- (z_\alpha)_j|$ for all $j\in\Lambda$. Hence, we have $ \lnorm{x^\dagger-T_\alpha(x^\dagger)}{\wei a}2 \leq \lnorm {x^\dagger-z_\alpha}{\wei a}2$.\\
It remains to find a bound on $ \lnorm {x^\dagger-z_\alpha} {\wei a}2$ 
in terms of $\|F(x^\dagger)-F(x_\alpha)\|_\Yspace$.\\ 
Let $\alpha>0$, $\beta:=2 L^2\alpha$ and $z_\alpha$ given by \eqref{eq:approximation_for_converse}. Then
\[ \frac{1}{2} \lnorm {x^\dagger- z_\alpha}{\wei a}2^2 + \alpha \lnorm {z_\alpha}{\wei r}1 \leq \frac{1}{2} \lnorm {x^\dagger- x_\beta}{\wei a}2^2 + \alpha \lnorm {x_\beta}{\wei r}1.\] 
Using \Cref{ass_operator} and subtracting $\alpha \lnorm {z_\alpha}{\wei r}1$ yield 
\begin{align} \label{eq:converse_eq2}
\frac{1}{2} \lnorm {x^\dagger- z_\alpha}{\wei a}2^2\leq \frac{L^2}{2}\|F(x^\dagger)-F(x_\beta)\|_\Yspace^2 + \alpha \left(\lnorm {x_\beta}{\wei r}1-\lnorm {z_\alpha}{\wei r}1 \right).
 \end{align}
Due to  \Cref{ass_domain} we have $z_\alpha\in D_F$.
As $x_\beta$ is a minimizer of \eqref{eq:Tik} we obtain 
\[ \beta \lnorm{x_\beta} {\wei r}1 \leq \frac{1}{2} \|F(x^\dagger)-F(x_\beta) \|_\Yspace^2 + \beta \lnorm{x_\beta}{\wei r}1 \leq \frac{1}{2} \|F(x^\dagger)-F(z_\alpha) \|_\Yspace^2 + \beta \lnorm{z_\alpha}{\wei r}1.\]
Using the other inequality in \Cref{ass_operator} and subtracting $\beta\lnorm {z_\alpha}{\wei r}1$ and dividing by $\beta$ we end up with 
\[ \lnorm{x_\beta} {\wei r}1 - \lnorm {z_\alpha}{\wei r}1 \leq \frac{L^2}{2\beta} \lnorm {x^\dagger-z_\alpha} {\wei a}2^2 = \frac{1}{4\alpha}  \lnorm {x^\dagger - z_\alpha} {\wei a}2^2.   \] 
We insert the last inequality into \eqref{eq:converse_eq2}, subtract $\frac{1}{4}  \lnorm {x^\dagger - z_\alpha} {\wei a}2^2$, multiply by $4$ and take the square root and get $\lnorm {x^\dagger-z_\alpha}{\wei a}2 \leq \sqrt{2} L \|F(x)-F(x_\beta)\|_\Yspace.$
Together with the first step, the hypothesis $(iii)$ and the definition of $\beta$ we achieve
\begin{align*}
\lnorm{x^\dagger-T_\alpha(x^\dagger)}{\wei a}2 \leq   \|F(x)-F(x_\beta)\|_\Yspace \leq (2 L^2)^\frac{3-t}{2}  C_3  \alpha^\frac{2-t}{2}. 
\end{align*}
Finally, \Cref{threshold_rates} yields $x\in k_t$ with $\|x^\dagger\|_{k_t}\leq c C_3^\frac{2}{t}$ with a constant $c$ that depends only on $t$ and $L$. 
\end{proof}

\section{Convergence analysis for $x^\dagger\notin \lspace {\wei r}1$}\label{sec:oversmoothing}
We turn to the oversmoothed setting where the unknown solution $x^\dagger$ does not admit a finite penalty value. 
An important ingredient of most variational convergence proofs of Tikhonov regularization is a comparison of the Tikhonov functional 
at the minimizer and at the exact solution. In the oversmoothing case such a comparison is obviously not useful. 
As a substitute, one may use a family of approximations of $x^\dagger$ at which the penalty functional is finite. See also  
\cite{HM:18} and \cite{HP:19} where this idea is used and the approximations are called auxiliary elements. Here we will use $T_{\alpha}(x^\dagger)$ for this purpose. 
We first show that the spaces $k_t$ can not only be characterized in terms of the approximation errors $\lnorm{(I-T_{\alpha})(\cdot)}{\wei \omega_p} p$ as in \Cref{threshold_rates}, 
but also in terms of $\lnorm {T_\alpha \cdot }{\wei r}1$:
\begin{lem}[bounds on $\lnorm {T_\alpha \cdot }{\wei r}1.$]\label{lem:threshold_bounds}
Let $t\in (1,2)$ and $x\in \mathbb{R}^\Lambda$. Then $x\in k_t$ if and only if $\eta(x):= \sup_{\alpha>0} \alpha^{t-1}\lnorm {T_\alpha(x)}{\wei r}1 <\infty$. \\
More precisely, we can bound 
\[ \eta(x)\leq 2 (1-2^{1-t})\inv \|x\|_{k_t}^t \quad\text{and }\|x\|_{k_t}\leq \eta(x)^\frac{1}{t}. \]
\end{lem} 
\begin{proof} 
As in the proof of \Cref{threshold_rates} we use a partitioning.  Assuming $x\in k_t$ we obtain 
\begin{align*}
\lnorm {T_\alpha(x)}{\wei r}1 & = \sum_{j\in\Lambda} \wei r_j |x_j| \mathds{1}_{ \{ \wei a_j^{-2} \wei r_j \alpha <|x_j| \} }   \\ 
& = \sum_{k=0}^\infty \sum_{j\in\Lambda} \wei r_j |x_j| \mathds{1}_{ \{ \wei a_j^{-2} \wei r_j 2^k \alpha < |x_j|\leq \wei a_j^{-2} \wei r_j 2^{k+1} \alpha \} }  \\ 
& \leq \alpha \sum_{k=0}^\infty  2^{k+1}\sum_{j\in\Lambda}\wei a_j^{-2} \wei r_j^2  \mathds{1}_{ \{ \wei a_j^{-2} \wei r_j 2^k \alpha < |x_j|\} }  \\ 
&\leq \|x\|_{k_t}^t\alpha^{1-t}  \sum_{k=0}^\infty 2^{k+1} 2^{-kt}  \\ 
&= 2 (1-2^{1-t})\inv \|x\|_{k_t}^t\alpha^{1-t}.
\end{align*}
Vice versa we estimate 
\begin{align*}
\sum_{j\in \Lambda}  \wei a_j^{-2}\wei r_j^{2}  \mathds{1}_{ \{ \wei a_j^{-2} \wei r_j\alpha< |x_j|\}} & \leq \alpha\inv\sum_{j\in \Lambda}  \wei r_j  |x_j| \mathds{1}_{ \{ \wei a_j^{-2} \wei r_j\alpha \leq |x_j|\}} \\ 
&= \alpha\inv \lnorm {T_\alpha(x)}{\wei r}1 \leq \eta(x) \alpha^{-t}.
\end{align*}
Hence, $\|x\|_{k_t}\leq \eta(x)^\frac{1}{t}.$
\end{proof} 

%A main ingredient is the following lemma evoking variational source conditions.
%\begin{lem} 
%Suppose \Cref{ass_operator} holds true. Let $p\in (0,2]$ and $x^\dagger\in k_p$. Then 
%\begin{align}\label{eq:replace_vsc}
%\lnorm {T_\alpha(x^\dagger)-x}{\wei r}1 +\lnorm {T_\alpha(x^\dagger)}{\wei r}1- \lnorm x{\wei r}1 \leq 2L \varrho_p(x^\dagger)^\frac{p}{2} \alpha^{-\frac{p}{2}} \|F(x^\dagger)-F(x)\|_\Yspace 
%\end{align}
%for all  $\alpha>0$ and  $x\in D\cap\lspace {\wei r}1$. 
%\end{lem} 
%\begin{proof}
%The coordinate projection $P_\alpha$ as defined in \Cref{lemma:bernstein} satisfies $P_\alpha x^\dagger= T_\alpha(x^\dagger).$ 
%\begin{align*}
%\lnorm {T_\alpha(x^\dagger)-x}{\wei r}1 +\lnorm {T_\alpha(x^\dagger)}{\wei r}1- \lnorm x{\wei r}1  & \leq \lnorm {P_\alpha x^\dagger -x}{\wei r}1 +\lnorm {P_\alpha x^\dagger}{\wei r}1- \lnorm x{\wei r}1  \\ 
%&= \lnorm {P_\alpha( x^\dagger -x
%)}{\wei r}1   +\lnorm {P_\alpha x^\dagger}{\wei r}1- \lnorm {P_\alpha x} {\wei r}1 \\ 
%&\leq 2 \lnorm {P_\alpha (x^\dagger -x) }{\wei r}1 \\
%&\leq 2 \varrho_p(x^\dagger)^\frac{p}{2} \alpha^{-\frac{p}{2}} \lnorm {x^\dagger-x} {\wei a}2 \\
%& \leq 2L \varrho_p(x^\dagger)^\frac{p}{2} \alpha^{-\frac{p}{2}} \|F(x^\dagger)-F(x)\|_\Yspace.  
%\end{align*}
%\end{proof} 
The following lemma provides a bound on the minimal value of the \linebreak Tikhonov functional. From this we deduce bounds on the distance between $T_\alpha(x^\dagger)$ and the minimizers of \eqref{eq:Tik} in  $\lnorm \cdot {\wei a}2$ and in  $\lnorm \cdot {\wei r}1.$ 
\begin{lem}[preparatory bounds]\label{lem:bounds}  Let $t\in (1,2)$, $\delta \geq 0$ and $\varrho>0$. 
Suppose \ref{ass_operator} and \ref{ass_domain} hold true. Assume $x^\dagger\in D_F$ with $\|x^\dagger\|_{k_t}\leq \varrho$ and $\gobs\in \Yspace$ with $\|\gobs-F(x^\dagger)\|_\Yspace \leq \delta.$ Then there exist constants $C_{t}$, $C_{a}$ and $C_{r}$ depending only on $t$ and $L$ such that
\begin{align}
\frac{1}{2} \|\gobs-F(\xh)\|_\Yspace^2  +\alpha \lnorm{\xh}{\wei r}1 & \leq    \delta^2 +  C_{t} \varrho^t \alpha^{2-t} \label{eq:bound_on_Tik}, \\ 
\lnorm  { T_\alpha(x^\dagger) - \xh }{\wei a}2^2 & \leq 8L^2 \delta^2 + C_{a} \varrho^t \alpha^{2-t}\quad\text{and}  \label{eq:bound_in_a2} \\ 
\lnorm { T_\alpha(x^\dagger) - \xh }{\wei r}1 &\leq \delta^2\alpha\inv+ C_{r}\varrho^t \alpha^{1-t}. \label{eq:bound_in_r1}
\end{align}
for all $\alpha>0$ and $\xh$ minimizers of \eqref{eq:Tik}.  
%\begin{thmlist}
%\item For all $c>1$ there exits $C_2>0$ depending only on $p$, $L$ and $c$ such that 
%\[  \|\gobs-F(\xh)\|_\Yspace^2  \leq c \delta^2 + C_2  \varrho_p(x^\dagger)^p \alpha^{2-p} \quad\text{for all }
%\alpha>0\] 
%\item 
%\end{thmlist}
\end{lem} 
\begin{proof}
Due to \Cref{ass_domain} we have $T_\alpha( x^\dagger)\in D$. Therefore, we may insert $T_\alpha( x^\dagger)$ into \eqref{eq:Tik} to start with
\begin{align}\label{eq:aux_element_inserted}
\frac{1}{2} \|\gobs-F(\xh)\|_\Yspace^2  +\alpha \lnorm{\xh}{\wei r}1 \leq  \frac{1}{2} \|\gobs -F(T_\alpha( x^\dagger))\|_\Yspace^2  +\alpha \lnorm{T_\alpha( x^\dagger)}{\wei r}1.
\end{align}
%We subtract $\alpha \lnorm{x_\alpha}{\wei r}1$ on both sides and use \eqref{eq:replace_vsc} with deleted first summand on the left hand side to wind up with 
%%\[ \frac{1}{2L^2} \lnorm {x^\dagger-x_\alpha} {\wei a}2 ^2 +\alpha \lnorm{x_\alpha}{\wei r}1 \leq  \frac{L^2}{2} \lnorm{ x^\dagger - P_\alpha x^\dagger } {\wei a}2 ^2  +\alpha \lnorm{ P_\alpha x^\dagger}{\wei r}1.\] 
%%We subtract $\lnorm{P_\alpha x_\alpha}{\wei r}1$ on both sides and use the reversed triangle inequality on the right hand side to wind up with 
%\begin{align}\label{eq:proof_of_over_rate}
%\frac{1}{2} \|\gobs-F(\xh)\|_\Yspace^2 \leq \frac{1}{2} \| \gobs - F(T_\alpha( x^\dagger))\|_\Yspace ^2  + 2L \varrho^\frac{p}{2} \alpha^\frac{2-p}{2} \|F(x^\dagger)-F(\xh)\|_\Yspace  . 
%\end{align}
\Cref{lem:threshold_bounds} provides the bound $ \alpha \lnorm{T_\alpha( x^\dagger)}{\wei r}1\leq C_1 \varrho^t \alpha^{2-t}$ for the second summand on the right hand side with a constant $C_1$ depending only on $t$. \\
In the following we will estimate the first summand on the right hand side. 
Let $\varepsilon>0$. 
By the second inequality in \Cref{ass_operator} and \Cref{threshold_rates} we obtain 
\begin{align}\label{eq:aux_prep}
\begin{aligned}
\frac{1}{2} \| \gobs - F(T_\alpha(x^\dagger))\|_\Yspace ^2 & \leq \| \gobs - F( x^\dagger)\|_\Yspace ^2 + \|F( x^\dagger) - F(T_\alpha(x^\dagger))\|_\Yspace ^2 \\ 
&\leq   \delta^2 +L^2 \lnorm{x^\dagger - T_\alpha(x^\dagger)} {\wei a}2 ^2\\ 
&\leq \delta^2  +  C_2  \varrho^t \alpha^{2-t}
\end{aligned}
\end{align}
with a constant $C_2$ depending on $L$ and $t$. Inserting into \eqref{eq:aux_element_inserted} yields \eqref{eq:bound_on_Tik} with $C_{t}:=C_1+C_2$. \\
%Choosing $\varepsilon=1$ we obtain 
%\[ \frac{1}{2} \| \gobs - F(T_\alpha(x^\dagger))\|_\Yspace ^2 \leq \delta^2 + 2C_2 \varrho^p \alpha^{2-p}\] 
%from the last inequality and 
We use \eqref{eq:aux_prep}, 
the first inequality in \Cref{ass_operator} and neglect the penalty term in \eqref{eq:bound_on_Tik} to estimate
\begin{align*}
\lnorm  { T_\alpha(x^\dagger) - \xh }{\wei a}2^2 &  \leq L^2 \| F(T_\alpha(x^\dagger)) - F(\xh) \|_\Yspace^2  \\  
& \leq 2L^2 \| \gobs - F(T_\alpha(x^\dagger))\|_\Yspace^2 + 2L^2 \|\gobs-F(\xh)\|_\Yspace^2 \\ 
&\leq 8L^2 \delta^2 + C_{a} \varrho^t \alpha^{2-t} 
\end{align*}
with $C_{a}:= 4L^2(C_2+C_{t})$.\\
\Cref{lem:threshold_bounds} provides the bound $ \lnorm {T_\alpha(x^\dagger)}{\wei r}1 \leq C_3 \varrho^t \alpha^{1-t}$ with $C_3$ depending only on $t.$
Neglecting the data fidelity term in \eqref{eq:bound_on_Tik}  yields
\begin{align}\label{eq:boundinsource}
\lnorm {T_\alpha(x^\dagger)-\xh}{\wei r}1 \leq  \lnorm {T_\alpha(x^\dagger)}{\wei r}1 +\lnorm {\xh} {\wei r}1 \leq \delta^2 \alpha\inv +C_r \varrho^t  \alpha^{1-t}
\end{align}  
with $C_r:=C_t+C_3.$
%Now assume $0<\varepsilon<1.$
%We use \Cref{lemma:bernstein}, the first inequality in \ref{ass_operator} and Young's inequality $2 ab\leq  \frac{4}{\varepsilon} \varepsilon\inv  a^2+ \frac{\varepsilon}{4} b^2$ to bound the second summand    
%\begin{align*}
%2L \varrho^\frac{p}{2} \alpha^\frac{2-p}{2} \|F(x^\dagger)-F(\xh)\|_\Yspace & \leq 4\varepsilon\inv L^2 \varrho^p \alpha^{2-p} +\frac{\varepsilon}{4} \|F(x^\dagger)-F(\xh)\|_\Yspace ^2 \\ 
%& \leq 4\varepsilon\inv L^2 \varrho^p \alpha^{2-p} +\frac{\varepsilon}{2} \delta^2 + \frac{\varepsilon}{2}  \|\gobs-F(\xh) \|_\Yspace^2
%\end{align*}
%Inserting the last two inequalities into \eqref{eq:proof_of_over_rate}, subtracting $ \frac{\varepsilon}{2}  \|\gobs-F(\xh)\|_\Yspace^2 $ and multiplying with $\frac{2}{1-\varepsilon}$ yields 
%\begin{align*}
% \|\gobs-F(\xh)\|_\Yspace^2  \leq \frac{1+2\varepsilon}{1-\varepsilon} \delta^2 + \frac{1+\varepsilon\inv}{1-\varepsilon} C_2  \varrho^p \alpha^{2-p}
%\end{align*}
%with $C_2=2 C_1 +8L^2.$
%The choice $\varepsilon:= \frac{c-1}{2+c}$ proves the claim.

\end{proof}
The next result is a converse type result for image space bounds with exact data. In particular, we see that Hölder type image space error bounds are determined by Hölder type bounds on the whole Tikhonov functional at the minimizers and vice versa. 
\begin{thm}[converse result for exact data]\label{thm:converse_oversmooth}
Suppose \Cref{ass_operator} and \ref{ass_domain} hold true. Let $t\in (1,2)$, $x^\dagger\in D_F$ and $(x_\alpha)_{\alpha>0}$ a choice of minimizers in \eqref{eq:Tik} with $\gobs=F(x^\dagger)$. The following statements are equivalent: 
\begin{enumerate}[(i)]
\item $x^\dagger\in k_t$. 
\item There exists a constant $C_2>0$ such that 
$\frac{1}{2}\|F(x)-F(x_\alpha)\|_{\Yspace}^2 +\alpha \lnorm{x_\alpha}{\wei r}1 \leq C_2 \alpha^{2-t}.$
\item There exists a constant $C_3$ such that $\|F(x)-F(x_\alpha)\|_{\Yspace}\leq C_3 \alpha^\frac{2-t}{2}$.
\end{enumerate}
More precisely, we can choose $C_2 = C_t \|x^\dagger\|_{k_t}^t$ with $C_t$ from \Cref{lem:bounds}, $C_3=\sqrt{2C_2}$ and bound $\|x^\dagger\|_{k_t}\leq c C_3^\frac{2}{t}$ with a constant $c$ that depends only on $t$ and $L$.
\end{thm}
\begin{proof}
\item[$(i)\Rightarrow (ii)$:] 
%Again the starting point is to  insert $T_\alpha(x^\dagger)$ into the Tikhonov functional: 
%\[ \frac{1}{2}\|F(x^\dagger)-F(x_\alpha)\|_\Yspace^2 + \alpha\lnorm{x_\alpha} {\wei r}1 \leq \frac{1}{2} \|F(x^\dagger)-F(T_\alpha(x^\dagger))\|_\Yspace^2 + \alpha \lnorm{T_\alpha(x^\dagger)}{\wei r}1. \]
Use \eqref{eq:bound_on_Tik} with $\delta=0$.
%\[ \frac{1}{2} \|F(x^\dagger)-F(T_\alpha(x^\dagger))\|_\Yspace^2+\alpha \lnorm{x_\alpha}{\wei r}1 \leq C_i \varrho_p(x^\dagger)^p \alpha^{2-p}.\]  
\item[$(ii)\Rightarrow (iii)$:] 
This implication follows immediately by neglecting the penalty term, multiplying by $2$ and taking the square root of the inequality in the hypothesis. 
\item[$(iii)\Rightarrow (i)$:] The same argument as in the proof of the implication (iii) $\Rightarrow$ (i) in \Cref{thm:converse} applies. 
%$ \lnorm{x^\dagger-T_\alpha(x^\dagger)}{\wei a}2 \leq   \|F(x)-F(x_\beta)\|_\Yspace \leq (2 L^2)^\frac{3-p}{2}$ for $\beta=2L^2 \alpha$. Hence, statement $(iii)$ together with \Cref{lem:threshold_bounds} yield $(i)$. 
\end{proof}

%\begin{lem}[error bounds in terms $\alpha$ and $\delta$] 
%Suppose \Cref{ass_operator} and \ref{ass_domain} hold true.
%Let $p\in (1,2)$, $e\in (p,2]$, $\varrho>0$ and $\delta\geq 0.$ Assume $x^\dagger\in D$ with $\varrho_p(x^\dagger)\leq \varrho$ and $\gobs\in \Yspace$ satisfies $\|\gobs-F(x^\dagger)\|_\Yspace\leq\delta$. There is a constant $C_b$ depending only $e, p$ and $L$ such that  
%\[ \] 
%\end{lem} 
%\begin{proof}
%Neglecting the data fidelity term in \eqref{eq:bound_on_Tik}  yields 
%\[ \lnorm {\xh}{\wei r}1 \leq \delta^2 \alpha\inv +   C_i \varrho^p \alpha^{1-p}\]  with $C_i$ from \Cref{lem:image_bounds}. \Cref{lem:threshold_bounds} provides the bound $ \lnorm {T_\alpha(x^\dagger)}{\wei r}1 \leq C_1 \varrho^p \alpha^{1-p}$ with $C_1$ depending only on $p.$ Hence
%\begin{align}\label{eq:boundinsource}
%\lnorm {T_\alpha(x^\dagger)-\xh}{\wei r}1 \leq  \lnorm {T_\alpha(x^\dagger)}{\wei r}1 +\lnorm {\xh} {\wei r}1 \leq \delta^2 \alpha\inv +(C_i+C_1) \varrho^p  \alpha^{1-p}
%\end{align}  
%
%\begin{align*}
%\lnorm {T_\alpha(x^\dagger)-\xh}e {\wei \omega_e} \leq 
%\end{align*}
%
%\[
%\lnorm{T_\alpha(x^\dagger)-\xh} {\wei a}2 \leq C_{ii}^\frac{1}{2} \varrho^\frac{p}{2} \alpha^\frac{2-p}{2}. 
%\] 
%
%\end{proof}
The following theorem shows that we obtain order optimal convergence rates on $k_t$ also in the case of oversmoothing  
(see \Cref{prop:optimal}). 
\begin{thm} [rates of convergence] \label{thm:rates_over}
Suppose Assumptions \ref{ass_operator} and \ref{ass_domain} hold true.
Let $t\in (1,2)$, $p\in (t,2]$ and $\varrho>0.$ Assume $x^\dagger\in D_F$ with $\|x^\dagger\|_{k_t}\leq \varrho$. 
\begin{enumerate}
\item (bias bound) Let $\alpha>0$.
For exact data $\gobs= F(x^\dagger)$ every minimizer $x_\alpha$ of \eqref{eq:Tik} satisfies 
\begin{align*}
\lnorm {x^\dagger-x_\alpha}{\wei \omega_p}p \leq C_b \varrho^\frac{t}{p} \alpha^\frac{p-t}{p}
\end{align*}  
with a constant $C_{b}$ depending only on $p, t$ and $L$.
\item  (rate with a-priori choice of $\alpha$)
Let $\delta>0$, $\gobs\in \Yspace$ satisfy $\|\gobs-F(x^\dagger)\|_\Yspace\leq\delta$ and $0<c_1<c_2$. If $\alpha$ is chosen such that \[ c_1 \varrho^\frac{t}{t-2} \delta^\frac{2}{2-t} \leq  \alpha \leq c_2 \varrho^\frac{t}{t-2} \delta^\frac{2}{2-t},\] then every minimizer $\hat{x}_\alpha$ of \eqref{eq:Tik} satisfies 
\begin{align*}
 \lnorm {\hat{x}_\alpha - x^\dagger} {\wei \omega_p}p \leq C_c  \varrho^\frac{t(2-p)}{p(2-t)}\delta^\frac{2(p-t)}{p(2-t)} 
\end{align*}  
with a constant $C_{c}$ depending only on $c_1, c_2, p, t$ and $L$.
\item  (rate with discrepancy principle) 
Let $\delta>0$ and $\gobs\in \Yspace$ satisfy $\|\gobs-F(x^\dagger)\|_\Yspace\leq\delta$ and $1<\tau_1\leq \tau_2$. If $\hat{x}_\alpha$ is a minimizer of \eqref{eq:Tik} with $\tau_1 \delta\leq \|F(\xh)-\gobs \|_\Yspace \leq \tau_2 \delta$, then 
\begin{align*}
\lnorm {\hat{x}_\alpha - x^\dagger} {\wei \omega_p}p\leq C_d  \varrho^\frac{t(2-p)}{p(2-t)}\delta^\frac{2(p-t)}{p(2-t)}.
\end{align*}
Here $C_d>0$ denotes a constant depending only on $\tau_1$, $\tau_2$, $p, t$ and $L$. 
\end{enumerate} 
\end{thm} 
\begin{proof}
\begin{enumerate}
\item 
By \Cref{prop:inter_scale} we have $\lnorm \cdot {\wei \omega_p}p \leq \lnorm \cdot {\wei a}2 ^\frac{2p-2}{p} \lnorm \cdot {\wei r}1^\frac{2-p}{p}.$ With this we interpolate between \eqref{eq:bound_in_a2} and \eqref{eq:bound_in_r1} with $\delta=0$ to obtain
\[\lnorm {T_\alpha(x^\dagger)-x_\alpha}{\wei \omega_p}p \leq  K_1 \varrho^\frac{t}{p} \alpha^\frac{p-t}{p}
\] 
with $K_1:=C_a^\frac{p-1}{p} C_r^\frac{2-p}{p}$. By \Cref{threshold_rates} there is a constant $K_2$ depending only on $p$ and $t$ such that  
\begin{align}\label{eq:thresh_e}
\lnorm {x^\dagger - T_\alpha(x^\dagger)}{\wei \omega_p}p \leq K_2\varrho^\frac{t}{p} \alpha^\frac{p-t}{p}.
\end{align}
 Hence
\begin{align*}
\lnorm {x^\dagger -x_\alpha}{\wei \omega_p}p & \leq \lnorm {x^\dagger - T_\alpha(x^\dagger)}{\wei \omega_p}p + \lnorm {T_\alpha(x^\dagger)-x_\alpha}{\wei \omega_p}p \\
& \leq (K_1+K_2) \varrho^\frac{t}{p} \alpha^\frac{p-t}{p}.
\end{align*} 
\item Inserting the parameter choice rule into \eqref{eq:bound_in_a2} and \eqref{eq:bound_in_r1} yields 
\begin{align*}
 \lnorm  { T_\alpha(x^\dagger) - \xh }{\wei a}2 & \leq (8L^2+C_a c_2^{2-t})^\frac{1}{2} \delta  \quad\text{and} \\ 
\lnorm { T_\alpha(x^\dagger) - \xh }{\wei r}1 &\leq (c_1\inv + C_r c_1^{1-t}) \varrho^\frac{t}{2-t}\delta^\frac{2(1-t)}{2-t}. 
\end{align*}
As above, we interpolate these two inequalities to obtain 
\[\lnorm {T_\alpha(x^\dagger)-\xh}{\wei \omega_p}p \leq K_3  \varrho^\frac{t(2-p)}{p(2-t)}\delta^\frac{2(p-t)}{p(2-t)}. 
\] 
with $K_3:= (8L^2+C_a c_2^{2-t})^\frac{p-1}{p} (c_1\inv + C_r c_1^{1-t})^\frac{2-p}{p}$. 
We insert the parameter choice into \eqref{eq:thresh_e} and get
$\lnorm {x^\dagger - T_\alpha(x^\dagger)}{\wei \omega_p}p \leq K_2 c_2^\frac{p-t}{p}\varrho^\frac{t(2-p)}{p(2-t)}\delta^\frac{2p-2t}{p(2-t)}.$ Applying the triangle inequality as in part $1$ yields the claim. 
\item Let $\varepsilon=\frac{\tau_1^2-1}{2}$. Then $\varepsilon>0$. By \Cref{threshold_rates} there exists a constant $K_4$ depending only on $t$ such that 
$\lnorm{x^\dagger-T_\beta(x^\dagger)}{\wei a}2^2 \leq K_4 \varrho^t \beta^{2-t}$ for all $\beta>0.$ We choose 
\[
\beta:=(\delta^2\varepsilon(1+\varepsilon\inv)\inv L^{-2} K_4\inv \varrho^{-t}  )^\frac{1}{2-t}.
\] 
Then 
\begin{align}\label{eq:bound_a2_dis}
\lnorm{x^\dagger-T_\beta(x^\dagger)}{\wei a}2^2 \leq \varepsilon(1+\varepsilon\inv)\inv L^{-2}\delta^2.
\end{align}
We make use of the elementary inequality 
$ (a+b)^2 \leq (1+\varepsilon)a^2 +(1 +\varepsilon\inv)b^2$ which is proven by expanding the square and 
applying Young's inequality on the mixed term. Together with the second inequality in \Cref{ass_operator} we estimate 
\begin{align*}
& \frac{1}{2} \| \gobs- F(T_\beta(x^\dagger)) \|_\Yspace^2 \\ & \leq \frac{1}{2} (1+\varepsilon) \|\gobs - F(x^\dagger)\|_\Yspace^2+ \frac{1}{2}(1+\varepsilon\inv) L^2 \lnorm{x^\dagger-T_\beta(x^\dagger)}{\wei a}2^2 \\ 
&\leq \frac{1}{2} (1+2\varepsilon) \delta^2 = \frac{1}{2} \tau_1^2 \delta^2. 
\end{align*}
By inserting $T_\beta(x^\dagger)$ into the Tikhonov functional we end up with 
\begin{align*}
\frac{1}{2} \tau_1^2 \delta^2 + \alpha \lnorm{\xh}{\wei r}1 & \leq \frac{1}{2} \|\gobs-F(\xh)\|_\Yspace^2 + \alpha \lnorm{\xh}{\wei r}1 \\  & \leq \frac{1}{2} \|\gobs-F(T_\beta(x^\dagger)) \|_\Yspace^2 +\alpha \lnorm{T_\beta(x^\dagger)}{\wei r}1\\& \leq  \frac{1}{2} \tau_1^2 \delta^2 + \alpha \lnorm{T_\beta(x^\dagger)}{\wei r}1.
\end{align*}
Hence, $\lnorm{\xh}{\wei r}1 \leq  \lnorm{T_\beta(x^\dagger)}{\wei r}1$. 
Together with \Cref{lem:threshold_bounds} we obtain the bound 
\[ \lnorm{T_\beta(x^\dagger)-\xh}{\wei r}1 \leq 2 \lnorm{T_\beta(x^\dagger)}{\wei r}1 \leq K_5 \varrho^\frac{t}{2-t}\delta^\frac{2-2t}{2-t} \] 
with a constant $K_5$ that depends only on $\tau$, $t$ and $L$.\\ 
Using \eqref{eq:bound_a2_dis} and the first inequality in \Cref{ass_operator} we estimate
\begin{align*}
& \lnorm{T_\beta(x^\dagger)-\xh}{\wei a}2 \\  & \leq \lnorm{x^\dagger - T_\beta(x^\dagger)}{\wei a}2+\lnorm{x^\dagger-\xh}{\wei a}2  \\ 
&\leq \lnorm{x^\dagger - T_\beta(x^\dagger)}{\wei a}2+L \|F(x^\dagger)-F(\xh)\|_\Yspace \\ 
&\leq  \lnorm{x^\dagger - T_\beta(x^\dagger)}{\wei a}2+L \|\gobs - F(x^\dagger)\|_\Yspace +L \|\gobs -F(\xh)\|_\Yspace \\
& \leq K_6 \delta
\end{align*} 
with $K_6= \varepsilon^\frac{1}{2}(1+\varepsilon\inv)^{-\frac{1}{2}}L^{-1}+L+L\tau_2.$ 
As above, interpolation yields 
\[\lnorm {T_\beta(x^\dagger)-\xh}{\wei \omega_p}p \leq K_7  \varrho^\frac{t(2-p)}{p(2-t)}\delta^\frac{2p-2t}{p(2-t)} 
\]
with $K_7:= K_6 ^\frac{2p-2}{p} K_5^\frac{2-p}{p}$. Finally, \Cref{threshold_rates} together with the choice of $\beta$ implies 
$\lnorm{x^\dagger-T_\beta(x^\dagger)} {\wei \omega_p}p \leq K_8 \varrho^\frac{t(2-p)}{p(2-t)}\delta^\frac{2p-2t}{p(2-t)}$ for a constant $K_8$ that depends only on $\tau$, $p,$ $t$ and $L$ and  we conclude 
\begin{align*}
 \lnorm{x^\dagger-\xh} {\wei \omega_p}p & \leq \lnorm{x^\dagger-T_\beta(x^\dagger)} {\wei \omega_p}p +\lnorm {T_\beta(x^\dagger)-\xh}{\wei \omega_p}p \\ & \leq (K_8+K_7) \varrho^\frac{t(2-p)}{p(2-t)}\delta^\frac{2p-2t}{p(2-t)}.
\end{align*}

\end{enumerate}
\end{proof}

\section{Wavelet Regularization with Besov Spaces Penalties}\label{sec:Besov}
\new{ In the sequel we apply our results developed in the general sequence space setting to obtain obtain convergence rates for wavelet regularization with a Besov $r,1,1$-norm penalty.}

%\subsection{Tikhonov Regularization with Besov $r,1,1$-Penalty}
Suppose Assumptions and \ref{ass:wavelet} and  \ref{ass_operator_Besov} and eqs.~\eqref{eqs:r_condition} hold true. 
Then  $F:=G\circ \wav$ satisfies \Cref{ass_operator} on $D_F:=\wav\inv(D_G) 
\subseteq \lspace {\wei a} 2= \bspace {-a}22$ 
as shown in \Cref{sec:ass_sequence}. \\
Recall that $\wei a _{(j,k)}= 2^{-ja}$ and $\wei r_{(j,k)}=2^{j(r-\frac{d}{2})}$.
Let $s\in [-a,\infty)$. With 
\begin{align}\label{eq:ts}
t_s:=\frac{2a+2r}{s+2a+r}
\end{align}
we obtain $\bspace s{t_s}{t_s}=\lspace {\wei \omega_{t_s}}{t_s}$ with equal norm for $\wei \omega_{t_s}$ given by \eqref{eq:omega}. For \new{$s\in (0,\infty)$ we have $t_s\in (0,1)$.}\\
\new{The following lemma defines and characterizes a function space $K_{t_s}$ as the counterpart of $k_{t_s}$ for $s>0$. 
As spaces $\bspace spq$ and $\Bspace spq$ with $p<1$ are involved let us first argue that within the scale $\bspace s{t_s}{t_s}$ for $s>0$ the extra condition $\sigma_{t_s}-\smax< s$ in \Cref{ass:wavelet} is always satisfied if we assume $a+r> \frac{d}{2}$. 
To this end let $0<s<\smax$. Then 
\begin{align*}
 \sigma_{t_s}= d\left( \frac{1}{t_s} -1 \right) = \frac{d(s-r)}{2a+2r}< s-r\leq s <\smax. 
\end{align*}
Hence, $\sigma_{t_s}- \smax < 0 <s$.
%Note that the condition $a+r> \frac{d}{2}$ is reasonable as it is equivalent to $\wei a_{(j,k)}\wei r_{(j,k)}\inv \rightarrow 0$, a part of \Cref{ass_operator}. 
}
\begin{lem}[Maximal approximation spaces $K_{t_s}$] 
\new{Let $a,s>0$ and suppose that \Cref{ass:wavelet} and eqs. \eqref{eq:r_condition} and \eqref{eq:r_nonneg} holds true. 
%Let $a>0$, $r\geq 0$ with $a+r>\frac{d}{2}$ and $s>0$. 
}
We define  \[ {K_{t_s}:=\mathcal{S}(k_{t_s})} \quad\text{ with } {\|f\|_{K_{t_s}}:=\|\mathcal{S}\inv x \|_{k_{t_s}}}\] with $t_s$ given by \eqref{eq:ts}. Let $s <u < s_\textrm{max}$. The space $K_{t_s}$ coincides 
with the real interpolation space
\begin{align}
\label{eq:Kt_interpol}
&K_{t_s}=(\Bspace {-a}22, \Bspace u {t_u}{t_u})_{\theta,\infty}, \qquad \theta = \frac{a+s}{u+a}.
\end{align}
with equivalent quasi-norms, and the following inclusions hold true with continuous embeddings:
\begin{align}\label{eq:Kt_embed}
&\Bspace s{t_s}{t_s} \subset K_{t_s} \subset \Bspace s{t_u}{\infty}.
\end{align}
Hence, \[ K_{t_s} \subset \bigcap_{t<t_s} \Bspace s t \infty.\]
\end{lem}
\begin{proof}
For $s<u<s_\textrm{max}$ we have $k_{t_s}=(\bspace {-a}22,\bspace u {t_u}{t_u} )_{\theta,\infty}$ with equivalent quasi-norms (see  \Cref{rem:weak}). By functor properties of real interpolation (see \cite[Thm.~3.1.2]{BL:76}) 
this translates to \eqref{eq:Kt_interpol}. \new {As discussed above, we use $a+r> \frac{d}{2}$ (see \eqref{eq:r_condition})
to see that $u\in (\sigma_{t_s}- \smax,s)$ such that $\mathcal{S}\colon \bspace  u {t_u}{t_u} \rightarrow \Bspace u {t_u}{t_u}$ is well defined an bijective.}
By \Cref{rem:markov} we have $\bspace s{t_s}{t_s}\subset k_{t_s}$ with continuous embedding, 
implying the first inclusion in \eqref{eq:Kt_embed}. Moreover, we have $t_u\leq \frac{2a+2r}{2a+r}\leq 2$. Hence, the continuous embeddings
$ \Bspace {-a}22 \subset \Bspace {-a}2 \infty \subset  \Bspace {-a}{t_u}\infty$ (see \cite[3.2.4(1), 3.3.1(9)]{Triebel2010}). Together with
\eqref{eq:Kt_interpol} and the interpolation result 
\[\Bspace s{t_u}{\infty} = (\Bspace{-a}{t_u}\infty,\Bspace u {t_u}{\infty})_{\theta,\infty}\] 
(see \cite[3.3.6 (9)]{Triebel2010}) we obtain the second inclusion in \eqref{eq:Kt_embed} using \cite[2.4.1 Rem.~4]{Triebel2010}. Finally, the last statement follows from $t_u\rightarrow t_s$ for $u\searrow s$ and again \cite[3.3.1(9)]{Triebel2010}.
\end{proof}

\begin{thm}[Convergence rates]\label{thm:Tik_besov}
Suppose Assumptions \ref{ass_operator_Besov} and \ref{ass:wavelet} hold true with $\frac{d}{2}-r<a<s_\textrm{max}$ and $\bspace r11 \cap \mathcal{S}\inv (D_G)\neq \emptyset$. Let $0<s<s_\textrm{max}$ with $s\neq r$, $\varrho>0$ and $\|\cdot\|_{L^p}$ denote the usual norm on $L^p(\Omega)$ for $1\leq p:=\frac{2a+2r}{2a+r}$. Assume $f^\dagger\in D_G$ with $\| f^\dagger\|_{K_{t_s}}\leq \varrho$. If $s<r$ assume that $D_F:=\wav\inv (D_G)$ satisfies \Cref{ass_domain}. Let $\delta>0$ and $ \gobs\in \Yspace$ satisfy $\|\gobs-F(f^\dagger)\|_\Yspace\leq \delta.$ 
\begin{enumerate}
\item  (rate with a-priori choice of $\alpha$)
Let $0<c_1<c_2$. If $\alpha$ is chosen such that \[ c_1  \varrho^{-\frac{a+r}{s+a}} \delta^\frac{s+2a+r}{s+a}\leq  \alpha \leq c_2 \varrho^{-\frac{a+r}{s+a}} \delta^\frac{s+2a+r}{s+a},\] then every $\hat{f}_\alpha$ given by \eqref{eq:Tik_besov} satisfies 
\begin{align*}
\left\|f^\dagger-\hat{f}_\alpha\right\|_{L^p}  \leq C_a  \varrho^\frac{a}{s+a}\delta^\frac{s}{s+a}.
\end{align*}  
\item  (rate with discrepancy principle) 
Let $1 < \tau_1 \leq \tau_2 $. If $\hat{f}_\alpha$ is given by   \eqref{eq:Tik_besov} with \[ \tau_1 \delta\leq \|F(\xh)-\gobs \|_\Yspace \leq \tau_2 \delta, \] then 
\begin{align*}
\left\|f^\dagger-\hat{f}_\alpha\right\|_{L^p}  \leq C_d  \varrho^\frac{a}{s+a}\delta^\frac{s}{s+a}.
\end{align*}
\end{enumerate}
Here $C_a$ and $C_{d}$ are constants independent of $\delta,$ $\varrho$ and $f^\dagger$.
\end{thm}
\begin{proof}
%Setting $a_{(j,k)}:= 2^{-ja}$ and $r_{(j,k)}:=2^{j(r-\frac{d}{2})}$ we obtain $\bspace r11=\lspace r1$ and $\bspace {-a}22=\lspace {\wei a}2$ with equal norms. \Cref{ass_operator_Besov} and \ref{ass:wavelet} imply that $F\circ \wav$   satisfies \Cref{ass_operator} on $D:=\wav\inv( D_F)$.\\
%Easy calculations prove $\bspace stt= \lspace {\omega_t} t$ and with equal norms and $\omega_t$ as defined in  \Cref{sec:weaksequencespace}. Likewise $\bspace 0pp= \lspace {\omega_p} p$ with equal norms.
%By \Cref{ass:wavelet} there is $x^\dagger\in \bspace stt$ with 
%$f^\dagger=\wav x^\dagger$. With $c_1$ the operator norm of $\wav\inv \colon \Bspace stt\rightarrow \bspace stt$ and in view of \Cref{rem:markov} we get 
%\[ \|x^\dagger\|_{k_t}\leq  \bn {x^\dagger} stt\leq c \Bn{f^\dagger} stt \leq c\varrho.\]
If $s>r$ (hence $t_s\in (0,1)$) we refer to \Cref{rem:error_bounds_intermediate}. If $s<r$ (hence $t\in (1,2)$) to \Cref{thm:rates_over} for the bound
\begin{align} \label{eq:boundproof}
\bn {x^\dagger - \xh} 0pp= \lnorm {x^\dagger - \xh} {\omega_p} p \leq C \varrho^{\frac{t_s}{p}\frac{2-p}{2-t_s}} \delta^{\frac{2}{p}\frac{p-t_s}{2-t_s}} = C \varrho^\frac{a}{s+a}\delta^\frac{s}{s+a}
\end{align}
for the a-priori choice
$ \alpha \sim \varrho^\frac{t_s}{t_s-2} \delta^\frac{2}{2-t_s}= \varrho^{-\frac{a+r}{s+a}} \delta^\frac{s+2a+r}{s+a} $ as well as for the discrepancy principle. 
%\item[2. case:] Suppose $s<r$. Then $t\in (1,2)$. We obtain the same bound as in the first case by 
%\item[3. case:] Suppose $s=r$. Then $t=1$ and $\lnorm{x^\dagger}r1 =\bn {x^\dagger}r11 \leq c\varrho.$ By we obtain a bound 
% \[ \bn {x^\dagger - \xh} 0pp= \lnorm {x^\dagger - \xh} {\omega_p} p \leq C \varrho^\frac{2-p}{p} \delta^\frac{2p-2}{p} = C  \varrho^\frac{a}{s+a}\delta^\frac{s}{s+a} \] 
%for both parameter choice rules. 
With \Cref{ass:wavelet} and by the well known embedding $\Bspace 0pp \subset L^p$ we obtain 
\[ \left\|f^\dagger-\hat{f}_\alpha\right\|_{L^p} \leq c_1 \Bn {f^\dagger-\hat{f}_\alpha}0pp \leq c_1 c_2  \bn {x^\dagger - \xh} 0pp.    \] 
Together with \eqref{eq:boundproof} this proves the result. 
\end{proof}
\begin{rem} 
In view of \Cref{rem:tequal1} we obtain the same results for the case $s=r$ by replacing $K_{t_s}$ by $\Bspace r11$. 
\end{rem}

\begin{thm}\label{thm:Besov_converse}
Let $r=0$. Suppose Assumptions \ref{ass_operator_Besov}, \ref{ass:wavelet} and \ref{ass_domain} hold true with $\smax >a>\frac{d}{2}$. Let $f^\dagger\in D_G\cap \Bspace 011$, $s>0$ and $(f_\alpha)_{\alpha>0}$ the minimizers of \eqref{eq:Tik_besov} for exact data $\gobs=F(f^\dagger)$. The following statements are equivalent: 
\begin{enumerate}[(i)] 
\item $f^\dagger \in K_{t_s}.$
\item There exists a constant $C_2>0$ such that $\Bn {f^\dagger-f_\alpha}011\leq C_2 \alpha^\frac{s}{s+2a}$ for all $\alpha>0$. 
\item There exists a constant $C_3>0$ such that $\|F(f^\dagger)-F(f_\alpha)\|_\Yspace\leq C_3 \alpha^\frac{s+a}{s+2a}$ for all $\alpha>0.$
\end{enumerate}   
More precisely, we can choose $C_2:= c \|f^\dagger\|_{K_t}^{t_s}$, $C_3:= c C_2^\frac{1}{2}$ and bound\linebreak ${\|f^\dagger\|_{K_t}\leq c C_3^\frac{2}{t_s}}$ with a constant $c>0$ that depends only on $L$ and $t$ and operator norms of $\mathcal{S}$ and $\mathcal{S}\inv$.  
\end{thm}
\begin{proof}
Statement $(i)$ is equivalent to $x^\dagger=\mathcal{S}\inv f^\dagger\in k_t$ and statement $(ii)$ is equivalent to a bound $\bn {x-x_\alpha}011\leq \tilde{C_2} \alpha^\frac{s}{s+2a}$. Hence, \Cref{thm:converse} yields the result. 
\end{proof}
\begin{examp}\label{ex:kink_jump}
\new{We consider functions $f^{\mathrm{jump}}, f^{\mathrm{kink}}:[0,1]\to \mathbb{R}$ which are $C^{\infty}$ everywhere with uniform bounds 
on all derivatives 
except at a finite number of points in $[0,1]$, and $f^{\mathrm{kink}}\in C^{0,1}([0,1])$. In other words, 
$f^{\mathrm{jump}}, f^{\mathrm{kink}}$ are piecewise smooth, $f^{\mathrm{jump}}$ has a finite number of jumps, and 
$f^{\mathrm{kink}}$ has a finite number of kinks. Then for $p\in (0,\infty)$, $q\in (0,\infty]$, and $s\in\mathbb{R}$ with $s>\sigma_p$ with $\sigma_p$ as in \Cref{ass:wavelet} we have
%which have an $l$th distributional derivative 
%of the form $f^{(l)}=\tilde{f}+\sum_{j=1}^N\gamma_j \delta_{x_j}$ with a smooth function $\tilde{f}$ and points $x_j\in(0,1)$. 
%In particular, for $l=1$ the function $f$ has jumps at the points $x_j$ and for $l=2$ it has kinks at these points. 
%Then by the definition of Besov spaces via the modulus of continuity (see, e.g., \cite[p.~3]{triebel:08}), we have 
\begin{align*}
f^{\mathrm{jump}}\in B^s_{p,q}((0,1)) \;\Leftrightarrow \; s<\tfrac{1}{p}, \qquad 
f^{\mathrm{kink}}\in B^s_{p,q}((0,1)) \; \Leftrightarrow \; s<1+\tfrac{1}{p} 
\end{align*}
if $q<\infty$ and 
\begin{align*}
f^{\mathrm{jump}}\in B^s_{p,\infty}((0,1)) \; \Leftrightarrow\;
s\leq \tfrac{1}{p},\qquad 
f^{\mathrm{kink}}\in B^s_{p,\infty}((0,1)) \; \Leftrightarrow\; 
s\leq 1+\tfrac{1}{p}.
\end{align*}
To see this, we can use the classical definition of Besov spaces in terms of the modulus of continuity $\|\Delta_h^m f\|_{L^p}$ 
where $(\Delta_hf)(x) := f(x+h)-f(x)$ and $\Delta_h^{m+1}f:= \Delta_h(\Delta_h^m f)$,  
see, e.g., \cite[eq.~(1.23)]{triebel:08}. Elementary computations show that 
 $\|\Delta_h^m f^{\mathrm{jump}}\|_{L^p}$ decays of the order $h^{1/p}$ as $h\searrow 0$ if $m\geq 1/p$, and 
 $\|\Delta_h^m f^{\mathrm{kink}}\|_{L^p}$ decays as $h^{1/p+1}$ if $m\geq 2/p$.}
\new{
Therefore, as $t_s<1$ describing the regularity of $f^{\mathrm{jump}}$ or $f^{\mathrm{kink}}$ in the scale $\Bspace s {t_s}{t_s} \subset K_{t_s}$ as in Theorems \ref{thm:Tik_besov} 
and \ref{thm:Besov_converse} allows for a larger 
value of $s$ and hence a faster convergence rate than describing the regularity of these functions in the Besov spaces $B^s_{1,\infty}$ as 
in \cite{HM:19}. 
In other words, the previous analysis in \cite{HM:19} provided only suboptimal rates of convergence for this important class 
of functions. This can also be observed in numerical simulations we provide below. }
\end{examp}

Note that the largest set on which a given rate of convergence is attained can be achieved by setting $r=0$ (i.e.\ no oversmoothing). 
This is in contrast to the Hilbert space case where oversmoothing allows to raise the finite qualification of Tikhonov regularization.  
On the other hand for larger $r$ convergence can be guaranteed in a stronger $L^p$-norm. 

%\todo[inline]{Eventuell optimale $L^p$-Raten für white noise bei $r>0$ (und somit $p>1$) diskutieren. Zusätzlich Lipschitz-Bed. an $F$ 
%nötig. Entscheidende Ungl.:
%\[
%\|g|B^{d/2}_{p,1}\| \leq \|g|B^0_{p,2}\|^{1-\frac{d}{2a}} \|g|B^a_{p,p}\|^{\frac{d}{2a}}
%\leq \|g|L^2\|^{1-\frac{d}{2a}} \|g|B^a_{p,p}\|^{\frac{d}{2a}}
%\]
%Erfordert aber doch recht viele neue Konzepte }

%\begin{corollary}
%In the setting of \Cref{thm:variant} assume $n$,$\alpha$ are chosen such that 
%\[ c_1 2^n \leq \varrho^\frac{1}{s+a} \delta^{-\frac{1}{s+a}} \leq c_2 2^n \quad\text{ and }\quad \alpha\leq c_3 \varrho^{-\frac{a+r}{s+a}} \delta^\frac{s+2a+r}{s+a}.\]
%Then there is a constant $\tilde{C}$ depending on $C$ in \Cref{thm:variant}, $c_1, c_2$ and $c_3$ such that 
% \[ \|f^\dagger- \hat{f}_\alpha\|_{L^2}\leq \tilde{C} \varrho^\frac{a}{a+s}\delta^\frac{s}{s+a}.\] 
%\end{corollary}
%\begin{proof}
%This follows immediately from plugging the parameter choices into the upper bound given in \Cref{thm:variant}. 
%\end{proof} 

\new{\section{Numerical results}\label{sec:num}
%We confirm our theoretical results in numerical experiments for the nonlinear parameter identification problem in \Cref{ex:identification}.  In one dimension it becomes the reconstruction of the coefficient $c$ in the elliptic 
%boundary value problem  
For our numerical simulations we consider the problem in \Cref{ex:identification} in the form 
\begin{align}\label{eq:bvp1d}
\begin{aligned}
&- u^{\prime\prime} + c u = f&&\mbox{in }(0,1),\\
&u(0)=u(1)=1.
\end{aligned}.
\end{align}
The forward operator in the function space setting is 
$G(c):=u$ for the fixed right hand side $f(\cdot)=\sin(4\pi \cdot )+2$.

%\old{
%For this problem \Cref{ass_operator_Besov} with $a=2$ can be proven based one the corresponding mapping property of the Fr\'echet derivative of $F$ (see \cite[Thm.~4.5]{HP:08}) and a \emph{range invariance condition} (see \cite[Ex.~4.2]{HNS:95}). We refer to \cite[Ex.~2.8, Lem.~2.9]{HM:19} for more details. } 
The true solution $c^\dagger$ is given by a piecewise smooth function with either finitely many jumps or kinks as discussed in  \Cref{ex:kink_jump}.

%Numerical simulations are carried out in Python.

To solve the  boundary value problem \eqref{eq:bvp1d} we used quadratic finite elements and an equidistant grid containing 
$127$ finite elements. The coefficient $c$ was sampled on an equidistant grid with $1024$ points. 
For the wavelet synthesis operator we used the code \texttt{PyWavelets} \cite{Lee2019} with Daubechies wavelet of order $7$. \\ 
The minimization problem in \eqref{eq:Tik_besov} was solved by the Gau{\ss}-Newton-type method $c_{k+1}= \wav x_{k+1}$, 
\begin{align*}
x_{k+1} \in \argmin_{x} \left[ \frac{1}{2} \| 	F^\prime[x_k](x-x_k)+ F( x_k)-u \|_\Yspace^2 +\alpha \bn {x-x_0} r11 \right] 
\end{align*} 
with a constant initial guess $c_0=1$. In each Gau{\ss}-Newton step 
these linearized minimization problems were solved with the \emph{Fast Iterative Shrinkage-Thresholding Algorithm 
(FISTA)}  proposed and analyzed by Beck \& Teboulle in \cite{beck2009fast}. We used the \emph{inertial parameter} as in \cite[Sec.~4]{chambolle2015convergence}. 
We did not impose a constraint on the size of $\bn{x-x_0} 022$, which is required by our theory if Assumption 3 does not hold true globally. However, the size of the domain of validity of this assumption is difficult to assess, and such a constraint is likely to be never active for a sufficiently good initial guess.

%For further numerical schemes to efficiently compute \eqref{eq:Tik} we refer to the overview paper \cite{CP:17}. 

The regularization parameter $\alpha$ was chosen by a sequential discrepancy principle with $\tau_1=1$ and $\tau_2=2$ on a grid $\alpha_j=2^{-j}\alpha_0$. To simulate worst case errors, we computed for each noise level $\delta$ reconstructions for several  data errors 
$u^{\delta}-G(c^\dagger)$, $\|u^{\delta}-G(c^\dagger)\|_{L^2}=\delta$, which were given by sin functions with different frequencies.}

\begin{figure}[ht]
\includegraphics[width=\textwidth]{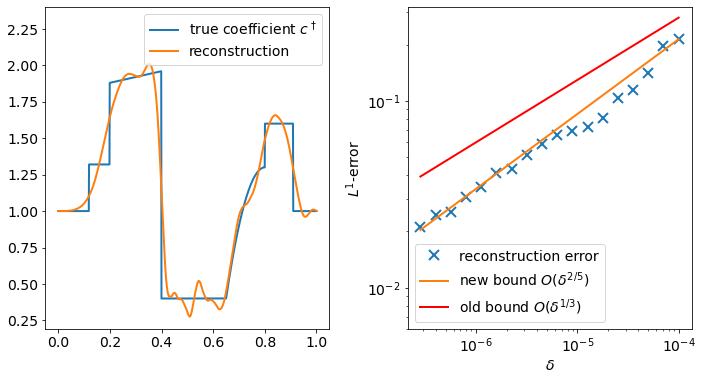} 
\caption{\label{fig:jump}\new{Left: true coefficient $c^\dagger$ with jumps in the boundary value problem \eqref{eq:bvp}  together with a typical reconstruction at noise level $\delta= 3.5\cdot 10^{-5}$. 
%\todo[inline]{Right: Also show $u^\dagger$ and $u^{\delta}$! Left: Legend entry 'reconstruction errors' 
%for crosses. Suggested order: reconstruction error, new bound, old bound. Second $y$-axis tick} 
Right: Reconstruction error using $\bspace 011$-penalization, the rate $\mathcal{O}(\delta^{2/5})$  
predicted by \cref{thm:Tik_besov} (see eq.~\eqref{eq:jump_rate}), and 
the rate $\mathcal{O}(\delta^{1/3})$ predicted by the previous analysis in \cite{HM:19}. 
}}
\end{figure}

\new{
For the piecewise smooth coefficient $c^\dagger$ with jumps shown on the left panel of Fig.~\ref{fig:jump}, 
\Cref{ex:kink_jump} yields 
\[
c^\dagger\in B^s_{t_s,t_s}((0,1))\subset K_{t_s} \quad \Leftrightarrow \quad s <\frac{1}{t_s} \quad \Leftrightarrow \quad  s<\frac{4}{3}. 
\]
Here $t_s=\frac{4}{s+4}$. Hence, \Cref{thm:Tik_besov} predicts the rate 
\begin{align}\label{eq:jump_rate}
\left\|c^\dagger-\widehat{c}_{\alpha}\right\|_{L^1} = \mathcal{O}(\delta^e)\qquad\mbox{ for all }e<\frac{2}{5}.
\end{align}
In contrast, the smoothness condition $c^\dagger \in B^s_{1,\infty} ((0,1))$ in our previous analysis in \cite{HM:19}, 
which was formulated in terms of Besov spaces with $p=1$ 
is only satisfied for smaller smoothness indices $s\leq 1$, and therefore, the convergence rate in \cite{HM:19} is only 
of the order $\left\|\widehat{c}_{\alpha}-c^\dagger\right\|_{L^1}= \mathcal{O}\left(\delta^{\frac{1}{3}}\right)$. 
Our numerical results displayed in the right panel of Fig.~\ref{fig:jump} show that this previous error bound is too pessimistic, 
and the observed convergence rate matches the rate \eqref{eq:jump_rate} predicted by our analysis. 
}

\begin{figure}[ht]
\centering{\includegraphics[width=\textwidth]{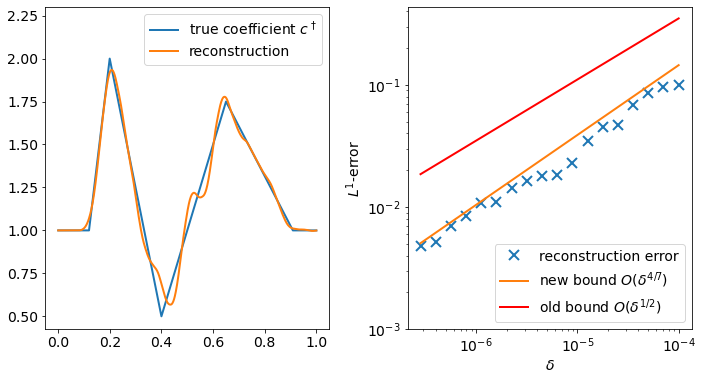}}
\caption{\label{fig:kink}\new{Left: true coefficient $c^\dagger$ with kinks in the boundary value problem \eqref{eq:bvp} 
together with a typical reconstruction at noise level $\delta= 3.5 \cdot 10^{-5}$. 
%\todo[inline]{As in Fig.~\ref{fig:jump}} 
Right: Reconstruction error using $\bspace 011$-penalization, the rate $\mathcal{O}(\delta^{4/7})$  
predicted by \cref{thm:Tik_besov} (see eq.~\eqref{eq:kink_rate}), and 
the rate $\mathcal{O}(\delta^{1/2})$ predicted by the previous analysis in \cite{HM:19}. 
}} 
\end{figure}

\new{Similarly, for the piecewise smooth coefficient $c^\dagger$ with kinks shown in the left panel of Fig.~\ref{fig:kink}, 
\Cref{ex:kink_jump} yields 
\begin{align*}
c^\dagger\in B^s_{t_s,t_s}((0,1))\subset K_{t_s} \quad \Leftrightarrow\quad
s <1+\frac{1}{t_s} \quad \Leftrightarrow \quad s<\frac{8}{3}
\end{align*}
with $t_s=\frac{4}{s+4}$.  Hence,  \Cref{thm:Tik_besov} predicts the rate 
\begin{align}\label{eq:kink_rate}
\left\|\widehat{c}_{\alpha}-c^\dagger\right\|_{L^1} = \mathcal{O}(\delta^e)\qquad \mbox{ for all }e<\frac{4}{7}
\end{align}
which matches with the results of our numerical simulations shown on the right panel of Fig.~\ref{fig:kink}. 
In contrast, the previous error bound $\left\|\widehat{c}_{\alpha}-c^\dagger\right\|_{L^1}=\mathcal{O}\left(\delta^\frac{1}{2}\right)$ 
in  \cite{HM:19} based on the regularity condition $c^\dagger \in B^2_{1,\infty} ((0,1))$ turns out to be suboptimal for this 
coefficient $c^\dagger$ even though it is minimax optimal in $B^2_{1,\infty}$-balls. 
}

\begin{figure}[ht]
\begin{center}
\includegraphics[width=\textwidth]{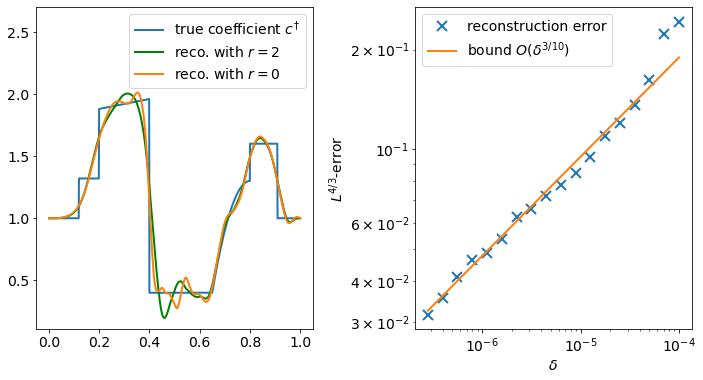}
\end{center}
\caption{\label{fig:reconstructions}\new{Left: true coefficient $c^\dagger$ with jumps in the boundary value problem \eqref{eq:bvp} 
together with reconstructions for $r=0$ and $r=2$ at noise level $\delta= 3.5\cdot 10^{-5}$ for the same data. 
Right: Reconstruction error using $\bspace 211$-penalization (oversmoothing) and the rate $\mathcal{O}(\delta^{3/10})$ predicted by \Cref{thm:Tik_besov} (see eq.~\eqref{eq:jump_over_rate}). This case is not covered by the theory in \cite{HM:19}.
}} 
\end{figure}

\new{Finally, for the same coefficient $c^\dagger$ with jumps as in Fig.~\ref{fig:jump}, reconstructions with $r=0$ and $r=2$ are compared in the left panel of Fig.~\ref{fig:reconstructions}. Visually, the reconstruction quality is similar for both reconstructions. 
For $r=2$ the penalization is oversmoothing, and \Cref{ex:kink_jump} yields 
\begin{align*}
c^\dagger\in B^s_{t_s,t_s}((0,1))\subset K_{t_s} \quad \Leftrightarrow\quad
s <\frac{1}{t_s} \quad \Leftrightarrow \quad s<\frac{6}{7}
\end{align*}
with $t_s=\frac{8}{s+6}$.  Hence,  \Cref{thm:Tik_besov} predicts the rate 
\begin{align}\label{eq:jump_over_rate}
\left\|\widehat{c}_{\alpha}-c^{\dagger}\right\|_{L^{4/3}} = \mathcal{O}(\delta^e)\qquad \mbox{ for all }e<\frac{3}{10},
\end{align}
which once again matches with the results of our numerical simulations shown on the right panel of Fig.~\ref{fig:reconstructions}. 
This case is not covered by the theory in \cite{HM:19}.}

\section{Conclusions}
We have derived a converse result for approximation rates of weighted $\ell^1$-regularization. 
Necessary and sufficient conditions for H\"older-type approximation rates are given by a scale of weak sequence spaces.
We also showed that $\ell^1$-penalization achieves the minimax-optimal convergence rates on bounded subsets of these 
weak sequence spaces, i.e.\ that no other method can uniformly perform better on these sets. 
However, converse results for noisy data, i.e.\ the question whether $\ell^1$-penalization achieves given convergence rates  
in terms of the noise level on even larger sets, remains open. Although it seems likely that the answer  
will be negative, a rigorous proof would probably require uniform
%A natural and interesting extension of this result would be an equivalent characterization of convergence rates for noisy data
%in terms of bounds in these weak sequence spaces. 
%We have only shown (using standard arguments) that the sufficient conditions for approximation rates also yield 
%convergence rates for noisy data. To prove the necessity of these conditions in this context, 
lower bounds on the maximal effect of data noise. 

A further  interesting extension concerns redundant frames. Note that \new{lacking injectivity} the composition of a forward operator in  function 
spaces with a synthesis operator of a redundant frame cannot meet 
the first inequality in Assumption \ref{ass_operator}. Therefore, the mapping properties of the forward 
operator in function space will have to be described in a different manner. \new{(See \cite[Sec.~ 6.2.]{AHR:13} for a related discussion.) }

We have also studied the important special case of penalization by wavelet Besov norms of type $B^r_{1,1}$. 
In this case the maximal spaces leading to H\"older-type approximation rates 
can be characterized as real interpolation spaces of
Besov spaces, but to the best of our knowledge they do not coincide with classical function spaces.
They are slightly larger than the Besov spaces $B^s_{t,t}$ with some $t\in (0,1)$, which in turn are considerably larger 
than the spaces $B^s_{1,\infty}$ used in previous results. 
Typical elements of the difference set $B^s_{t,t}\setminus B^s_{1,\infty}$ 
are piecewise smooth functions with local singularities. 
Since such functions can be well approximated by functions with sparse wavelet expansions, good performance 
of $\ell^1$-wavelet penalization is intuitively expected. Our results confirm and quantify this intuition.

\appendix
\section{Appendix}
For a sequence $(\wei \omega_j)_{j\in J}$ of positive real numbers, we write $\wei \omega_j\rightarrow 0$ if for every $\varepsilon>0$ the set $\{ j\in \Lambda \colon \wei \omega_j > \varepsilon \}$ is finite. 
\begin{prop}[embeddings] \label{app:embed}
Let $1\leq p \leq q <\infty$ and $s= (s_j)_{j\in \Lambda}$, $r= (\wei r_j)_{j\in \Lambda}$ sequences of positive reals. 
\begin{enumerate}[(i)] 
\item There is a continuous embedding $\lspace r {p} \subset \lspace s {q}$ iff $s_j \wei r_j\inv$ is bounded. 
\item There is a compact embedding $\lspace r {p} \subset \lspace s {q}$ iff $s_j \wei r_j\inv \rightarrow 0$. 
\end{enumerate}
\end{prop}
\begin{proof}
\item[(i)] If there is such a continuous embedding, then there exists a constant $C>0$ such that \linebreak ${\lnorm \cdot s {q} \leq C \lnorm \cdot r {p}.}$ Inserting unit sequences $e_j:=(\delta_{jk})_{k\in \Lambda}$ yields $s_j \wei r_j\inv \leq C$.\\ 
For the other implication we assume that there exists a constant $C>0$ such that $s_j \wei r_j\inv \leq C$ for all $j\in\Lambda$. Let $x\in \lspace r{p}$ with $ \lnorm x r p=1$. Then $s_j |x_j|\leq C \wei r_j |x_j| \leq C  \lnorm x r p $ implies 
\[\lnorm x s {q}^q =\sum_{j\in \Lambda} s_j^{q} |x_j|^q \leq (C \lnorm x r p)^{q-p} \sum_{j\in \Lambda} s_j^{p} |x_j|^{p} \leq C^q\lnorm x r p^{q-p} \sum_{j\in \Lambda} \wei r_j^p |x_j|^{p}= C^q \lnorm x r p^{q}.\]
Taking the $q$-th root shows $\lnorm \cdot s {q} \leq C \lnorm \cdot r {p}.$
\item[(ii)]
Suppose $s_j \wei r_j\inv \rightarrow 0$ is false. Then there exists some $\varepsilon$ and a sequence of indices $(j_k)_{k\in\mathbb{N}}$ such that $s_{j_k} r_{j_k}\inv\geq \varepsilon$ for all $k\in \mathbb{N}.$ The sequence given by $x_k=r_{j_k}\inv e_{j_k}$ is bounded in $\lspace rp$. But $\lnorm {x_k-x_m} s q \geq 2^\frac{1}{q} \varepsilon$ for $k\neq m$ shows that it does not contain a convergent subsequence in $\lspace sq$.\\
To prove the other direction we assume $s_j \wei r_j\inv \rightarrow 0$. Then $s_j \wei r_j\inv$ is bounded and by part (i) there is a continuous embedding $I\colon  \lspace rp \rightarrow \lspace sq$. We define $\Lambda_n=\{ j\in\Lambda \colon s_j \wei r_j\inv> \frac{1}{n}\}$. As $\Lambda_n$ is finite the coordinate projection $P_n\colon \lspace rp \rightarrow \lspace sq$ given by $(P_n x)_j= x_j$ if $j\in \Lambda_n$ and $(P_n x)_j= 0$ else is compact.    
As $ s_j \wei r_j\inv\leq \frac{1}{n}$ for all $j\in \Lambda\setminus\Lambda_n$ part (i) yields 
\[ \lnorm {(I-P_n)x} sq \leq \frac{1}{n} \lnorm {(I-P_n)x} rp \leq \frac{1}{n} \lnorm {x} rp \quad\text{for all } x\in \lspace rp.\]
Hence, $\|I-P_n\|\leq \frac{1}{n}$. Therefore, $I=\lim_n P_n$ is compact. 
\end{proof}
Financial support by Deutsche Forschungsgemeinschaft (DFG, German Science Foundation) 
through grant RTG 2088 is gratefully acknowledged.  %    amsalpha, or (for "historical" overv

%\begin{acknowledgements}
%\end{acknowledgements}


\bibliographystyle{spmpsci}     % mathematics and physical sciences
\bibliography{lit}   % name your BibTeX data base

\end{document}